\documentclass[final,1p,times]{elsarticle}
 \usepackage{amsfonts}

\usepackage{mathrsfs, amsmath,amssymb,bm}
\usepackage{graphicx}
\usepackage{subfigure}
\usepackage{caption}
\usepackage{overpic}
\usepackage{epsfig}
\usepackage{xcolor}
\usepackage{indentfirst}
\usepackage{fancyhdr}
\usepackage{comment}
\usepackage{makecell, rotating}
\usepackage{algorithm} 
\usepackage{float}
\usepackage{diagbox}
\usepackage{xkeyval}
\usepackage{todonotes}
\presetkeys{todonotes}{inline}{} 

\usepackage[draft]{changes}
\usepackage{bm}
\usepackage{float}
\usepackage{cite}
\usepackage{multirow}
\usepackage{booktabs}
\usepackage{xkeyval}
\usepackage{mathtools}
\usepackage{tikz}
\usetikzlibrary{shapes.geometric, arrows}
\showboxbreadth\maxdimen\showboxdepth\maxdimen

\newcommand{\be}{\bm{e}}

\newcommand{\bmm}{\bm{m}}

\newcommand{\bq}{\bm{q}}
\newcommand{\bu}{\bm{u}}
\newcommand{\bx}{\bm{x}}

\newcommand{\bF}{\boldsymbol{F}}

\newcommand{\bM}{\mathcal{M}}
\newcommand{\bV}{\mathcal{V}}
\newcommand{\bC}{\mathcal{C}}

\newcommand{\bK}{\mathcal{K}}

\newcommand{\calA}{\mathcal{A}}

\newcommand{\bbR}{\mathbb{R}}

\newcommand{\bbM}{\bm{M}}
\newcommand{\bA}{\bm{A}}

\usepackage[draft]{changes}
\usepackage{lipsum}
\definechangesauthor[name={Kai Jiang}, color=red]{JK}

\begin{document}

\title{An adaptive high-order surface finite element method for polymeric
 self-consistent field theory on general curved surfaces}
\author{Kai Jiang, Xin Wang, Jianggang Liu, Huayi Wei\corref{cor1}}
\cortext[cor1]{Corresponding author. Email: weihuayi@xtu.edu.cn}

\address{School of Mathematics and Computational Science,       
Hunan Key Laboratory for Computation and Simulation in Science and Engineering,
Xiangtan University, Xiangtan, Hunan, China, 411105.}
\date{\today}

\begin{abstract}
In this paper, we develop an adaptive high-order surface finite element method
(FEM) incorporating spectral deferred correction method for
chain contour discretization to solve polymeric self-consistent field equations
on general curved surfaces. The high-order surface FEM is obtained by
the high-order surface geometrical approximation and the high-order function space
approximation. Numerical results demonstrate that the
precision order of these methods is consistent with theoretical prediction.
In order to describe the sharp interface in the strongly
segregated system more accurately, an adaptive FEM equipped with a new
\textit{Log} marking strategy is proposed. Compared with the traditional
strategy, the \textit{Log} marking strategy can not only label the elements that need to
be refined or coarsened, but also give the refined or coarsened times, which can
make full use of the information of a posterior error estimator and improve the
efficiency of the adaptive algorithm.  To demonstrate the power of our approach,
we investigate the self-assembled patterns of diblock copolymers on several
distinct curved surfaces.  Numerical results illustrate the efficiency of the
proposed method, especially for strongly segregated systems with economical
discretization nodes. 
\end{abstract}

\maketitle
\section{Introduction}
\label{sec:intrd}

In recent years, microphase separation of block copolymers under various types
of geometrical confinements (including Euclidean and manifold confinements) has
attracted tremendous attention due to their industrial
applications\,\citep{stewart2011, segalman2005}. Geometric constraints
drastically affect the formation of ordered structures under which some
traditional ordered phases of block copolymers are rearranged to form novel
patterns\,\citep{stewart2011, segalman2005, wu2004confinement,
  xiang2005confinement, yu2006confinement, deng2016chiral}.  Theories play an
  important role in understanding and predicting the phase behavior of block
  copolymers under geometrical confinements\,\citep{zhang2014sphere,
  li2006self,chantawansri2007, wei2019finite, fredrickson2006equilibrium}. Among
  these theories, self-consistent field theory
  (SCFT)\,\citep{fredrickson2006equilibrium} is one of the most powerful tools in
  studying the self-assembly behaviors of inhomogeneous polymers and related
  soft-matter systems.

There have been a few studies done on the surface SCFT calculations,
versus numerous works on computing bulk structures based on
SCFT.
Chantawansri \textit{et al.}\,\citep{chantawansri2007} and Vorselaars \textit{et
al.}\,\citep{Vorselaars2011} used the spherical harmonic method to numerically
simulate the SCFT model confined on the spherical surface and in the spherical
shell, respectively. The global spherical harmonic
method has spectral accuracy, however, it can not extend to general curved surfaces. 
For computing the SCFT model on general surface,  Li \textit{et
al.}\,\citep{li2006self, li2014self} proposed a method similar to
the finite volume method.  However, there has been no theory (numerical result) to guarantee (rigorously
demonstrate) the computational precision.
Besides that, Li's method can not be applied to strongly segregated systems when
interaction parameter $\chi N > 25$ for diblock copolymer melt.
Meanwhile, even for relatively weak interaction systems $\chi N <25$, Li's
method still requires a large number contour discretization points (more than
$1000$) to reduce the free-energy discrepancy about $5\times 10^{-5}$.
Precisely computing strongly segregated systems is still a challenge in the
SCFT computation, especially for general surface confinement. 
In this work, we are devoted to developing efficient high-order numerical
methods for polymer SCFT on general curved surfaces.

In the past several decades, many numerical approaches have been developed to
address surface problems, including the level-set method\,\citep{osher2001}, the close point method \citep{macdonald2009}, and the surface
FEM\,\citep{Dziuk1988, demlow2007,Dziuk2013, Bonito2020}. 
In this paper, we focus on the surface FEM.
Dziuk\,\citep{Dziuk1988} firstly
proposed a linear FEM to solve the Laplace-Beltrami equation on arbitrary
surfaces.  Demlow and Dziuk\,\citep{demlow2007} presented an adaptive linear
surface FEM, and then Demlow\,\citep{Demlow2009} generalized the surface FEM theory
to a high-order case. Wei \textit{et al.}\citep{wei2010} generalized the superconvergence
results and several gradient recovery approaches of linear FEM from flat spaces
to general curved surfaces for the Laplace-Beltrami equation with mildly
structured triangular meshes.  Bonito and Demlow\,\citep{Bonito2019} gave the new
posteriori error estimates with sharper geometric error estimators for surface
FEMs. After about 25 years of development, the surface FEM has been applied
to a wide range of scientific problems, see recent review
papers\,\citep{Dziuk2013, Bonito2020} and references therein.  

In our previous work\,\citep{wei2019finite}, we proposed a linear surface FEM to
study the microphase separation of block copolymers on general curved surfaces.
However, in SCFT calculations, using the linear surface FEM to achieve a
relatively high numerical precision may result in heavy computational
complexity. Meanwhile, in strong segregation regime, self-assembled structures
have two-scale spatial distribution: sharp interfaces and damped internal densities,
making the uniform mesh method inefficient.  Therefore it
is necessary to improve this numerical method. 
The main contributions of this work include:
\begin{itemize}
	\item We present an adaptive high-order surface FEM for polymeric SCFT.  
		It is high order both in space (surface discretization
		and function space approximation) and in time through incorporating
		spectral deferred correction method for chain contour integration. 
		A concrete way to construct arbitrary order surface FEM is also given.
	\item We propose a novel and efficient marking strategy in the adaptive method.
		This new marking approach does not only denote which mesh element needs
		to be changed, but also provides the times of refinement or coarseness.
		Compared with existing marking strategies, the new approach can efficiently
		make full use of the information of a posterior error estimator.
	\item Our developed method is attractive for strongly segregated systems.  
		Numerical results for strongly segregated systems demonstrate that our proposed
		approach can achieve the prescribed precision with an economical
		computational cost both in space and time directions.
\end{itemize}

The rest of this paper is organized as follows. In Sec.\,\ref{sec:model}, we
introduce a self-consistent field model on general curved surfaces of diblock
copolymer melt as an example to demonstrate our method.  In
Sec.\,\ref{sec:method}, we present the concrete construction procedure of
adaptive high-order surface FEM in detail and apply it to solve the propagator
equation. In Sec.\,\ref{sec:contour}, we incorporate the spectral deferred
correction scheme for the contour discretization into the adaptive high-order FEM. 
In Sec.\,\ref{sec:rslt}, we examine the efficiency of the proposed methods
through sufficient numerical examples.  It should be noted that the proposed
method is also suitable for other polymeric systems. 

\section{Surface self-consistent field theory}
\label{sec:model}

In this section, we present the SCFT for an incompressible
AB diblock copolymer melt on a generally curved surface using the standard Gaussian chain model. 
We consider $n$ AB diblock copolymers confined on a curved surface $\bM$ whose
measure is $|\bM|$. The volume fraction of A block 
is $f$ and that of B block is $1-f$, the total degree of
polymerization of a diblock copolymer is $N$. The field-based Hamiltonian for the
incompressible diblock copolymer melt is\,\citep{fredrickson2006equilibrium}
\begin{equation}\label{eq:hamiltonian}
    H = \frac{1}{|\bM|}\int_{\bM}\left\{-\omega_+(\bx) + \frac{\omega_-^2(\bx)}{\chi
N}\right\}\mathrm{d}\bx -\log Q[\omega_+(\bx), \omega_-(\bx)],
\end{equation}
where $\chi$ is the Flory-Huggins parameter to measure the interaction
between segments A and B. The term $\omega_+(\bx)$ is the fluctuating pressure
field, and $\omega_-(\bx)$ is the exchange chemical potential field. The
pressure field enforces the local incompressibility, while the exchange
chemical potential field is conjugate to the density operators. $Q$ is the single
chain partition function.

First-order variations of the Hamiltonian with respect to the fields $\omega_\pm(\bx)$
lead to the following mean-filed equations,
\begin{align}
	\label{eq:scft1}
\frac{\delta H}{\delta \omega_+(\bx)} & =\varphi_A(\bx)+\varphi_B(\bx)-1=0,
\\
	\label{eq:scft2}
\frac{\delta H}{\delta \omega_-(\bx)}
&= \frac{2\omega_-(\bx)}{\chi N}-[\varphi_A(\bx)-\varphi_B(\bx)]=0.
\end{align}
$\varphi_A(\bx)$ and $\varphi_B(\bx)$ are the monomer densities of blocks A and
B, respectively.
\begin{equation}\label{eq:phiA}
\varphi_A(\bx)=\frac{1}{Q}\int_0^fq(\bx,t)q^{\dagger}(\bx,t) \,\mathrm{d}t,
\end{equation}
\begin{equation}\label{eq:phiB}
\varphi_B(\bx)=\frac{1}{Q}\int_f^1q(\bx,t)q^{\dagger}(\bx,t) \, \mathrm{d}t.
\end{equation}
The single chain partition function $Q$ can be calculated by
\begin{equation}\label{eq:Q}
    Q  = \frac{1}{|\bM|}\int_{\bM} q(\bx,1)\, \mathrm{d}\bx.
\end{equation}
The forward propagator $q(\bx,t)$ denotes the probability distribution of 
the chain in contour node $t$ at surface position $\bx$. The variable
$t$ is used to parameterize each copolymer chain such that $t=0$ represents the
tail of the A block and $t=f$ is the junction between the A and B blocks,
$t=1$ is the end of the B block. From the continuous Gaussian chain
model\,\citep{matsen2002standard}, $q(\bx,t)$ satisfies the PDE
\begin{equation}\label{eq:forward}
\begin{aligned}
    \frac{\partial}{\partial t}q(\bx,t)&=[\Delta_{\bM}-w(\bx,t)]q(\bx,t),\\
q(\bx,0)&=1,\\
w(\bx,t) &= \left\{   
\begin{array}{rl}
\omega_+(\bx)-\omega_-(\bx), &  0\leq t \leq f,\\
\omega_+(\bx)+\omega_-(\bx), &  f\leq t \leq 1.
\end{array}
\right.
\end{aligned}
\end{equation}
$\Delta_{\bM}$ is the Laplace-Beltrami operator which reads as follows
$$
\Delta_{\bM} v = div_{\bM} \cdot (\nabla_{\bM} v),
$$
where $\bM$ is a two-dimensional, compact, and $\bC^2$-hypersurface in $\bbR^3$.
$\nabla_{\bM}$ and $div_{\bM}$ are the tangential
gradient and divergence operators, respectively, their definitions can be found in\,\citep{Dziuk1988, demlow2007}. The Sobolev space on surface $\bM$ is 
$$
H^p(\bM) = \{ v\in L^2(\bM)\ |\ D^{\alpha} v \in L^2(\bM),\ |\alpha| \leq p\}.
$$

The backward propagator $q^{\dagger}(\bx,t)$ represents the probability
distribution from $t=1$ to $t=0$ satisfying
\begin{equation}\label{eq:backward}
\begin{aligned}
\frac{\partial}{\partial t}q^{\dagger}(\bx,t)
&=[\Delta_{\bM}-w^{\dagger}(\bx,t)]q^{\dagger}(\bx,t),\\
q^{\dagger}(\bx,0)&=1,\\
w^{\dagger}(\bx,t) &= \left\{   
\begin{array}{rl}
\omega_+(\bx)+\omega_-(\bx), &  0\leq t \leq 1-f,\\
\omega_+(\bx)-\omega_-(\bx), &  1-f\leq t \leq 1.
\end{array}
\right.
\end{aligned}
\end{equation}
For closed surfaces, no boundary condition is needed.
While for open surfaces, Eqns.\,\eqref{eq:forward} and \eqref{eq:backward}
require boundary conditions to be well-posed.
In this work, we use the homogeneous Neumann boundary condition for open surfaces. Certainly,
other appropriate boundary conditions can be considered in the following proposed algorithm framework.

\begin{figure}[!hbpt]
	\begin{center}	
		\tikzstyle{startstop} = [rectangle, rounded corners, minimum width=3cm,
        minimum height=1cm, text width=3cm, text centered, draw=black,
      fill=blue!20]
		\tikzstyle{io} = [rectangle, minimum width=4cm,minimum
        height=1cm,text width=6cm,text centered,draw=black,fill=green!20]
		\tikzstyle{ioo} = [rectangle, rounded  corners,minimum width=3cm,minimum
        height=1cm,text centered,draw=black,fill=red!30]
		\tikzstyle{process1} = [rectangle,  minimum
        width=3cm,minimum height=1cm,text width=7cm, text
      centered,draw=black,fill=green!20]
		\tikzstyle{process2} = [rectangle, minimum
        width=3cm,minimum height=1cm,text width=5.5cm, text
      centered,draw=black,fill=green!20]
		\tikzstyle{decision} = [trapezium, trapezium left angle=70, trapezium right angle=110,  rounded corners, minimum width=3cm,minimum
      height=1cm,text centered,text width=5cm, draw=black,fill=orange!30]
        \tikzstyle{arrow} = [thick=50cm,->,>=stealth]
		\begin{tikzpicture}[node distance=1.7cm]
		minimum	\node (start) [startstop] {Give initial fields $\omega_+(\bx)$ and $\omega_-(\bx)$};
		\node (input1) [io,below of=start] {Calculate propagators $q(\bx,t)$ and
        $q^\dag(\bx, t)$ on a general curved surface $\bM$};
		\node (process1) [process1,below of=input1,yshift=-0.1cm] {Obtain $Q$,
        density operators $\varphi_A(\bx)$ and $\varphi_B(\bx)$, and calculate the Hamiltonian $H$};
		\node (process2) [process2,below of=process1,yshift=-0.1cm] {Update
        fields $w_+(\bx)$ and $w_-(\bx)$};
        \node (decision) [decision, below of =process2,yshift=-0.1cm] {Hamilton
        discrepancy is lower than a prescribed criterion};
		\node (out1) [ioo,below of=decision] {Stop};
        \coordinate (point1) at (5cm, -7.1cm);
        \draw [arrow] (start) -- (input1);
		\draw [arrow] (input1) -- (process1);
		\draw [arrow] (process1) -- (process2);
      \draw (decision) -- node [above] {N} (point1);
      \draw [arrow] (point1) |- (input1);
        \draw [arrow] (decision) --node [right] {Y} (out1);
		\draw [arrow] (process2) -- (decision);
		\end{tikzpicture}
	\end{center}
	\caption{Flowchart of SCFT iteration.}
	\label{fig:scftiter}
\end{figure}
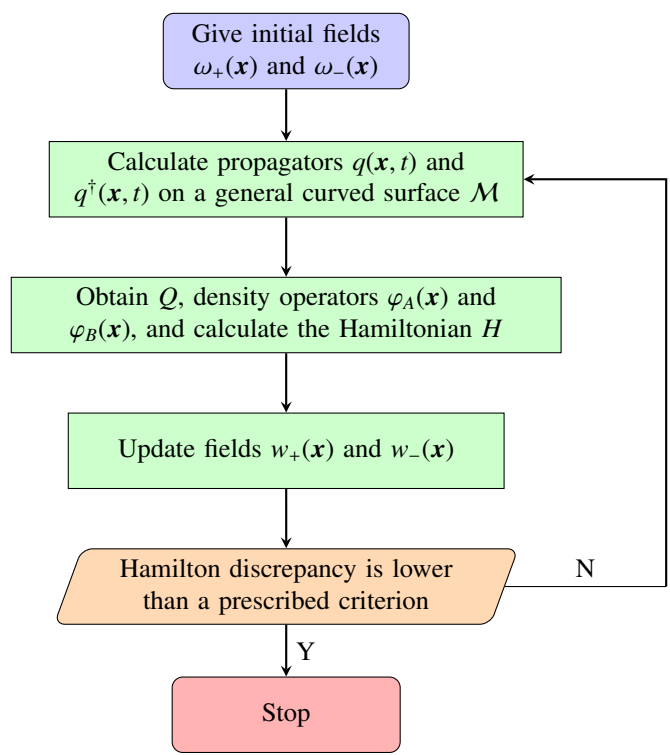

 Finding equilibrium states of SCFT corresponds to solutions of Euler-Lagrange
equation $\delta H/\delta \omega_{\pm}(\bx)=0$. Note that field functions
$\omega_{\pm}(\bx)$ are related to the density functions $\varphi_A(\bx)$, $\varphi_B(\bx)$ which satisfy integral
equations \eqref{eq:phiA} and \eqref{eq:phiB}. The integrands in
Eqns.\,\eqref{eq:phiA} and \eqref{eq:phiB} are related to PDEs \eqref{eq:forward} and
\eqref{eq:backward}. Thus it is a nonlocal problem. Usually, iteration methods
are designed to update field functions $\omega_\pm(\bx)$. 
The standard SCFT iteration procedure is shown in the Fig.\,\ref{fig:scftiter}.

It is important to note that the iteration approach for updating the fields
$\omega_+(\bx)$ and $\omega_-(\bx)$ could be determined by the SCFT mathematical properties.
The Hamiltonian $H$ of AB diblock copolymers can reach its local minima along
the exchange chemical field $\omega_-(\bx)$ and achieve the maxima along the pressure
field $\omega_+(\bx)$\,\citep{fredrickson2006equilibrium}. For multicomponent
polymer systems, an extended analysis of the SCFT model can be found in
Ref.\,\citep{jiang2015analytic}.
Several iteration methods have been proposed to find the saddle point, such as the
semi-implicit Seidel method\,\citep{fredrickson2004iteration}, which is based on the asymptotic expansion and global
Fourier transformation. Similarly, the semi-implicit scheme using global
spherical harmonic
transformation can be extended to the spherical surface
problems\,\citep{chantawansri2007}. 
However, it may be impossible to obtain a global
basis on general curved surfaces. Correspondingly, the semi-implicit scheme
could not be obtained. Thus, we choose the alternative direction explicit
Euler method to update fields, i.e., 
\begin{equation}\label{iteration}
	\begin{aligned}
		\omega_+^{k+1}(\bx) &= \omega_+^k(\bx) + \lambda_+
		[\varphi_A^k(\bx)+\varphi_B^k(\bx)-1 ],
		\\
		\omega_-^{k+1}(\bx) &= \omega_-^k(\bx) - \lambda_-
        \left(\frac{2 \omega_-^k(\bx)}{\chi N} - [\varphi^k_A(\bx) - \varphi^k_B(\bx)]
   \right),
	\end{aligned}
\end{equation}
where $\lambda_{\pm}$ are the iteration step size.

\section{Adaptive high-order surface FEM}
\label{sec:method}



For the strongly segregated polymer systems, the self-assembled patterns have sharp
interfaces and gradually changed internal structures. A uniform mesh method
gives rise to a large amount of computational complexity.  This
section presents
an adaptive high-order surface FEM to solve the propagators on general
curved surfaces, obtaining better accuracy results with lower computation cost.

Consider a curved surface $\bM$ with boundary $\bM_g$ on which \eqref{eq:forward} and
\eqref{eq:backward} are well-defined. Let $\bV=H^1(\bM)$ be the trial and test
function spaces. The variational formulation is stated as follows:
find $q\in \bV$ such that for $\forall v\in\bV$ 
\begin{equation}
	\label{eq:variation}
	\left(\frac{\partial}{\partial t} q, v\right)_{\bM} = -(\nabla_{\bM} q,
	\nabla_{\bM} v)_{\bM} - (wq, v)_{\bM}, ~~\text{ for all }~ v\in \bV,
\end{equation}
where $(\cdot,\cdot)_{\bM}$ represents the $L^2$ inner product on curved
surface $\bM$. 

\subsection{The construction of high-order surface FEM space}
Next we present the construction of the high-order surface FEM
space. We first introduce the multi-index vector $\bmm:=(m_0, m_1, m_2)$ with  
\begin{equation*}
	m_i \geq 0, ~ i=0, 1, 2, \text{ and } \sum_{i=0}^2 m_i=p.
\end{equation*}
The number of all possible of $\bmm$ is 
\begin{equation*}
	n_p := \begin{pmatrix}
		p+2 \\ 2 
	\end{pmatrix}.
\end{equation*}

\begin{figure}[!hbpt]
	\centering	
    \includegraphics[width=1.0\linewidth]{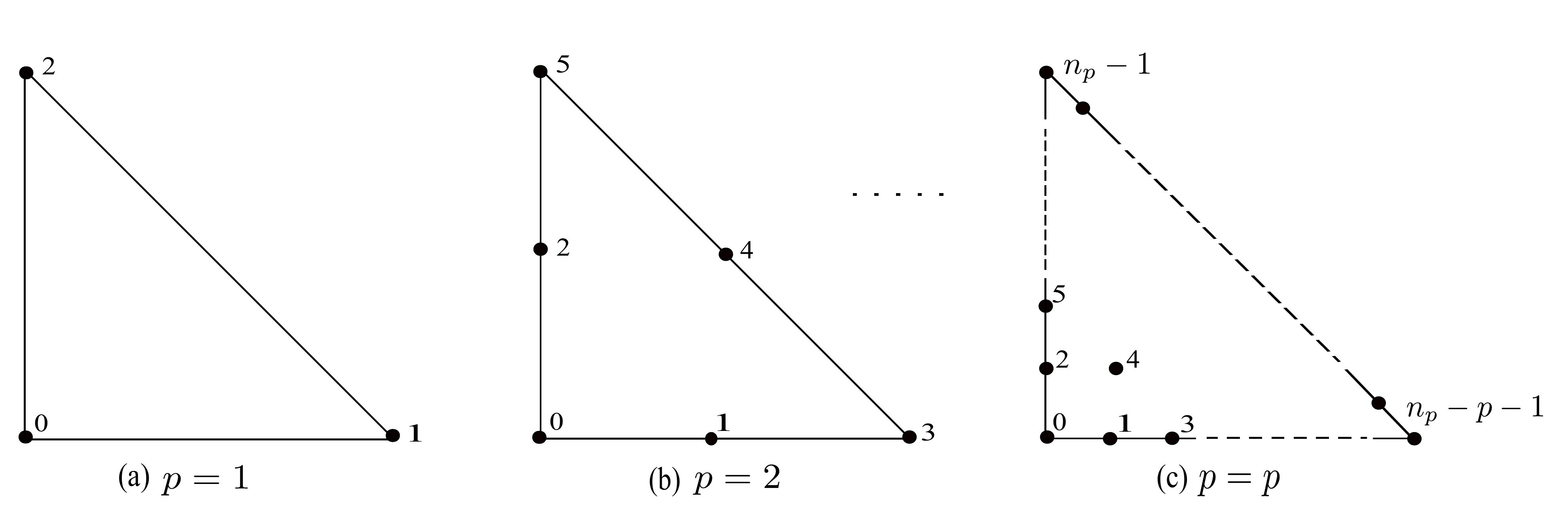}
    \caption{The spatial sorting rules for the  points $\bu_{\alpha}$  on the
    reference triangle $\hat{\tau}$.}
\label{fig:tau}
\end{figure}

Let $\alpha$ be a one-dimensional index starting from 0 to $n_p-1$ of
$\bmm$, see Tab.\,\ref{tb:num} for the numbering rule.
\begin{table}[H]
	\centering
	\caption{The numbering rule of multi-index $\bmm_\alpha$.}\label{tb:num}
	\begin{tabular}{| l | c | c | c|}
	\hline
	$\alpha$ & \multicolumn{3}{c|}{$\bmm_\alpha$} \\
	\hline
	$0$ & $p$   & $0$ & $0$ \\\hline
	$1$ & $p-1$ & $1$ & $0$ \\\hline
	$2$ & $p-1$ & $0$ & $1$ \\\hline
	$3$ & $p-2$ & $2$ & $0$ \\\hline
	$4$ & $p-2$ & $1$ & $1$ \\\hline
	$\vdots$ & $\vdots$ & $\vdots$ & $\vdots$ \\\hline
	$n_p-1$ & $0$ & $0$ & $p$ \\
	\hline
	\end{tabular}
\end{table}
For each multi-index $\bmm_\alpha=(m_0, m_1, m_2)$, one can define a point
$\bu_\alpha=(\xi,\eta)=(m_1/p, m_2/p)$ on the
reference triangle $\hat{\tau}$ (an isosceles right triangle with the right side
length of $1$), and the corresponding barycentric coordinate $\bm\lambda =
(\lambda_0, \lambda_1, \lambda_2)$, where
\begin{equation}
  \lambda_0 = \frac{m_0}{p} = 1-\xi-\eta, \quad \lambda_1 = \frac{m_1}{p} = \xi,
  \quad \lambda_2 = \frac{m_2}{p} = \eta.
\end{equation}
Using the numbering rule, the spatial sorting of
$\bu_\alpha$ on $\hat{\tau}$ follows the rules shown in Fig.\,\ref{fig:tau}.
Then for each multi-index $\bmm_\alpha$, one can construct the shape function of
degree $p$ on the reference triangle $\hat{\tau}$ as follows
\begin{equation}
  \phi_{\alpha}(\bu) = \frac{1}{\bmm_\alpha!}\prod_{i=0}^{2}\prod_{j =
  0}^{m_i - 1} (p\lambda_i - j), 
	\label{eq:phi0}
\end{equation}
with   
\begin{align*}
	\bmm_\alpha! = m_0!m_1! m_2! \,, ~~ \text{ and }\prod_{j=0}^{-1}(p\lambda_i -
	j) = 1,\quad i=0, 1,2.
\end{align*}
It is easy to verify that the shape function $\phi_\alpha$ defined above
satisfies the interpolation property
\begin{equation}
	\phi_\alpha(\bu_\beta) = 
	\begin{cases}
		1, & \alpha = \beta,\\
		0, & \alpha \ne \beta,
	\end{cases}
	\quad \text{ with }\alpha, \beta = 0, 1, \cdots, n_p - 1.
	\label{eq:phi1}
\end{equation}
Obviously, $\{\phi_\alpha\}_{\alpha=0}^{n_p-1}$ are linearly independent.
Notice that a formula similar to \eqref{eq:phi0} can be found
in\,\citep{nguyen2010}. Here we write it in a slightly different form and clearly
show out its coefficients. Obviously, this formula can be easily extended to any
dimensional geometric simplex and one can find such implementation in
FEALPy\,\citep{fealpy2021}.

\begin{figure}[!hbpt] \centering
  \includegraphics[width=0.3\linewidth]{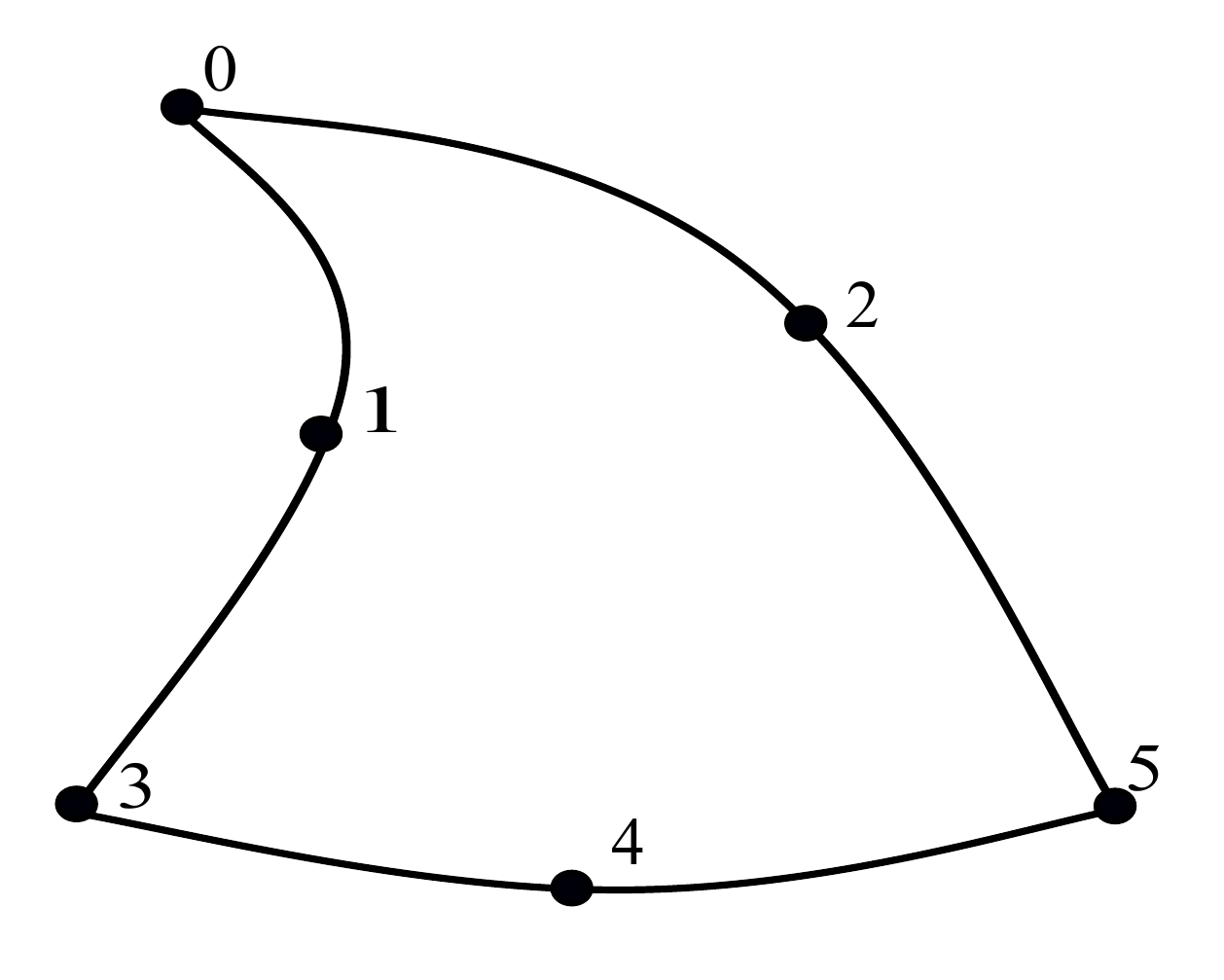} 
  \caption{The quadric triangle surface $\tau^2_h$.} \label{fig:tpmesh} 
\end{figure}

Let $\bM_h^p = \{\tau_h^p\}$ be a union of a set of triangle surfaces which is
a continuous piecewise polynomial approximation of degree $p$ of surface $\bM$.
Each triangle surface $\tau_h^p$ is uniquely determined by a set of
interpolation points $\{\bx_\alpha\}_{\alpha=0}^{n_p -1} \subset \bM$, which
follows the same numbering rule as $\{\bu_\alpha\}_{\alpha=0}^{n_p-1}$ on
reference element $\hat{\tau}$, see Fig \ref{fig:tpmesh} for the case of $p=2$.
Then we can give a one-to-one mapping $\bx: \hat{\tau} \mapsto \tau_h^p$ as follows
\begin{equation*} 
  \bx(\bu) = \sum_{\alpha=0}^{n_p-1}
  \bx_\alpha\phi_\alpha(\bu),\,\forall \bu\in \hat{\tau}, 
\end{equation*} 
and obviously,
\begin{equation*} 
  \bx_\alpha = \bx(\bu_\alpha), \alpha=0, 1, \cdots, n_p-1.
\end{equation*}

\begin{figure}[!hbpt] \centering
  \includegraphics[width=0.6\linewidth]{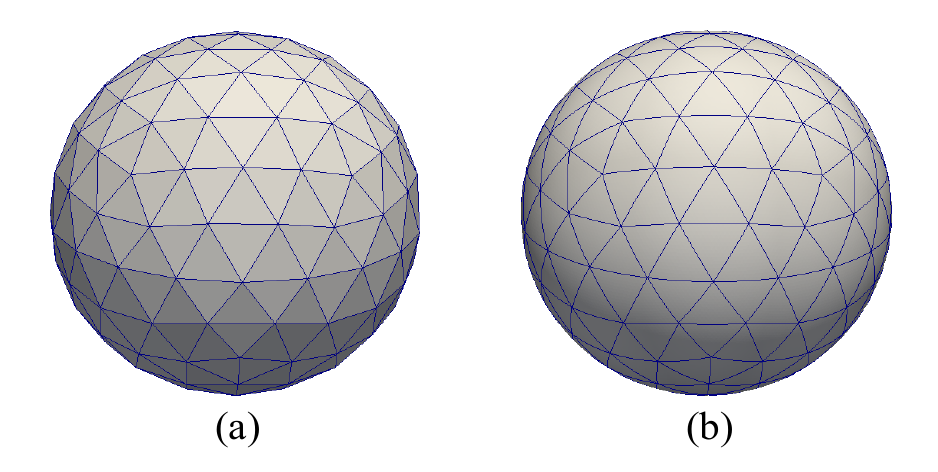} 
  \caption{The piecewise linear surface $\bM_h^1$ (a) and  quadric surface $\bM_h^2$ (b)
triangular mesh on a sphere.} \label{fig:mpmesh} 
\end{figure}

In Fig.\,\ref{fig:mpmesh}, we present a piecewise linear approximation surface
$\bM_h^1$ and a piecewise quadric surface $\bM_h^2$ of a sphere.

  By the one-to-one mapping $\bx: \hat{\tau}
\mapsto \tau_h^p$, we can define the local basis function on $\tau_h^p$
\begin{equation*}
  \varphi_\alpha(\bx) := \phi_\alpha(\bu),\quad \forall \bx(\bu) \in \tau_h^p \text{ with
  }
  \bu\in\hat\tau. 
\end{equation*}

Similarly, $\varphi_\alpha$ satisfies the interpolation property on
interpolation points $\{\bx_\alpha\}_{\alpha=0}^{n_p -1}$ of $\tau_h^p$
\begin{equation*}
	\varphi_\alpha(\bx_\beta) = 
	\begin{cases}
		1, & \alpha = \beta,\\
		0, & \alpha \ne \beta,
	\end{cases}
	\quad \text{ with }\alpha, \beta = 0, 1, \cdots, n_p - 1.
\end{equation*}
Furthermore, we know that the adjacent elements on high-order surface
$\bM_h^p$ will reuse the interpolation points
on their shared nodes or edges. Suppose $\bM_h^p$ has $NN$ nodes, $NE$ edges and $NC$
elements, the number of interpolation points on $\bM_h^p$ is  
\begin{equation*}
  N_{DoF}:= NN + (p-1)NE + \frac{(p-1)(p-2)}{2}NC.
\end{equation*}
For each interpolation point, one can construct a continuous piecewise global basis
function on $\bM_h^p$, which take value 1 at this interpolation
point and 0 at others. On the element the interpolation point is located, the 
global basis function is exactly the local basis function corresponding to the
interpolation point on that element. For simplicity, we still use $\varphi_i$ to represent the
$i$-th global basis function.

Finally, we introduce a continuous piecewise polynomial space of degree $p$ defined on
$\bM_h^p$   
\begin{equation*} 
  \bV_{h}^p := \text{span}\{\varphi_i\}_{i=0}^{N_{DoF}-1}.
\end{equation*} 
Notice
that, the above definition is equivalent to the definition in the literature
\citep{Dziuk2013}(Page 315, Intermediate remark).

The error of the surface FEM contains two parts: the geometric error arising
from the approximation of $\bM$ by the $\bM^p_h$ and the function space
approximation error coming from the approximation of an infinite-dimensional
function space by a finite-dimensional space. As Demlow
analyzed\,\citep{Demlow2009, libuyang2021}, when employing finite element space of degree $p$
on the polynomial surface approximation of degree $k$, one can have
\begin{align}
  \|q-q^l_h\|_{L^2(\bM)} \leq C h^{p+1} \|q\|_{H^{p+1}(\bM)} +
Ch^{k+1}\|q\|_{H^1(\bM)}, 
\end{align} where $q$ is the exact solution,
$q^l_h$ denotes the lift of numerical solution $q_h$ from $\bM^p_h$ to
$\bM$, $C$ depends on geometric properties of
$\bM$. In this work, we set $p=k$ to ensure that the approximation errors of the
geometric and function space have the same accuracy.

\subsection{Surface FEM discretization}
By the finite-dimensional surface FEM
space $\bV_h^p$, the descrete weak formulation of  
\eqref{eq:variation} is: find
$$
q_h(\bx,t)=\sum_{i=0}^{N_{DoF}-1} q_i(t) \varphi_i(\bx)  \in \bV_h^p
$$
such that
\begin{equation}\label{eq:discretization}
	\left(\frac{\partial}{\partial t} q_h, v_h\right)_{\bM_h^p} 
	= -(\nabla_{\bM_h^p}q_h, \nabla_{\bM_h^p} v_h)_{\bM_h^p} - (\omega_hq_h,
    v_h)_{\bM_h^p},
~~\text{ for all } v_h\in \bV_h^p.
\end{equation}
Let $v_h = \varphi_j$, $j=0, 1, \cdots, N_{DoF}-1$, we have the matrix form of
\eqref{eq:discretization}
\begin{equation}\label{eq:matrix}
	\bbM\frac{\partial}{\partial t}\bq(t) = - (\bA + \bF)\bq(t),\\
\end{equation}
where 
\begin{equation*}
	\bq(t) = \left(q_0(t), q_1(t), \cdots, q_{N_{DoF}-1}(t)\right)^T,
\end{equation*} 
and 
\begin{equation*}
	\bbM_{i,j} = (\varphi_i, \varphi_j),\quad \bA_{i, j} = (\nabla_{\bM_h^p}\varphi_i,
	\nabla_{\bM_h^p}\varphi_j), \quad \bF_{i, j} = (\omega_h\varphi_i, \varphi_j).
\end{equation*}

\subsection{Adaptive surface FEM}
\label{subsec:afem}

\begin{figure}[!htbp]
	\centering
	\includegraphics[width=1.0\linewidth]{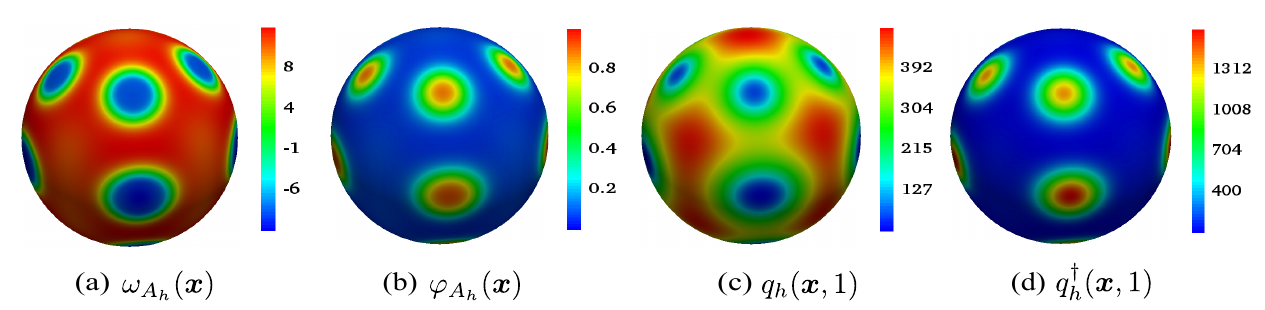}
	\caption{\label{fig:estimator} The equilibrium distributions of
      $\omega_{A_h}(\bx)$, $\varphi_{A_h}(\bx)$, $q_h(\bx,1)$ and $q_h^{\dagger}(\bx,1)$ on a
    spherical surface with a radius 3.56 $R_g$ when $\chi N =25, f=0.2$.}
\end{figure}

In this subsection, we present the adaptive mesh method to obtain a
high-precision numerical solution with less computational complexity.
The adaptive method can automatically rearrange mesh grids according to the
error distribution over each mesh element. 
The adaptive procedure used in the SCFT calculation contains
\begin{description}
	\item[Step 1] Solve the SCFT model and obtain the numerical solution on a
		given mesh.
	\item[Step 2] Calculate a posterior error estimator on each element from current
		numerical results.
	\item[Step 3] Mark mesh elements according to the error estimator.
	\item[Step 4] Refine or coarsen the marked elements to update the mesh.
	\item[Step 5] Go to Step 1 until the desired error is satisfied.
\end{description}
Define the error estimator of the indicator function $u_h$ on each element $\tau$ as
$e_{\tau}$,
\begin{equation}
	e_{\tau} = \|G_h u_h\|_{\tau},
	\label{eq:afem:error}
\end{equation}
where $\|\cdot\|_{\tau}$ is denoted as the $L^2$ norm on $\tau$. $G_h u_h$ is
the harmonic average operator\,\citep{huang2012}
\begin{equation}
	G_h u_h(\bx_i) := \frac{1}{\sum_{j=1}^{n}1/|
    \tau_j|}\sum_{j=1}^{n} \frac{1}{|\tau_j|}\nabla_{\bM_h^p} u_h \Big|_{\tau_j}(\bx_i).
\end{equation}
$n$ is the number of elements $\tau_j$ with $\bx_i$ as a vertex.
An appropriate indicator function is of importance in the adaptive mesh method.  
To obtain the appropriate indicator function in such complicated SCFT calculations,
we observe the distribution of different spatial functions of
the equilibrium state, including field function $w_{A_h}(\bx)$, density function
$\varphi_{A_h}(\bx)$, and propagators $q_h(\bx,t)$, $q_h^{\dagger}(\bx,t)$, such as at
the last contour point $q_h(\bx,1)$, $q_h^{\dagger}(\bx,1)$. 
As an example, Fig.\,\ref{fig:estimator} gives the corresponding results of
these spatial functions when $\chi N =25, f=0.2$. 
From these results, one can find that 
their spatial distributions are similar, however, $q_h(\bx,1)$ and
$q_h^{\dagger}(\bx,1)$ vary more dramatically and have sharper interfaces.
It can be expected that the error of numerical solutions decreases as  
the numerical errors of $q_h(\bx,1)$ and $q_h^{\dagger}(\bx,1)$ decrease.
Therefore using $q_h(\bx, 1)$ and $q_h^{\dagger}(\bx,1)$ as the indicator function
to obtain the posterior error estimator $e_{\tau}$ is a proper choice.
Concretely, $e_{\tau}$ can be chosen as

\begin{equation}                                                                
    e_{\tau} = \frac{\|G_h q_h(\bx,1)\|_{\tau}+\|G_h q_h^{\dagger}(\bx,1)\|_{\tau}}{2}.
    \label{eq:afem:error1}                                                       
\end{equation}

Besides the error estimator $e_{\tau}$, the mesh marking strategy 
is also crucial. Classic marking strategies, including the
maximum criterion\,\citep{jarausch1986} and the $L^2$
criterion\,\citep{dorfler1996}, usually refine
or coarsen marked elements one time in each adaptive process.
It may make less use of the information of the posterior error estimator. 
To improve it, we propose a new marking strategy, called the \textit{Log} criterion,
as follows
\begin{equation}
  \label{eq:log}
	n_{\tau} = \left[\log_2\frac{e_{\tau}}{\theta \bar{e}}\right],
\end{equation} 
where $\theta$ is a positive constant, $\bar{e}$ is the mean value of
estimators $e_{\tau}$ on all cells, and $[\cdot]$ is the nearest integer function.
$n_{\tau} = 0$, $n_{\tau} > 0$ and $n_{\tau} < 0$ represent that element $\tau$
is unchanged, refined $n_{\tau}$ times, and coarsened $|n_{\tau}|$ times, respectively.
This new marking strategy does not only denote which mesh element 
$\tau$ needs to be changed, but also provides the times of refinement or coarseness. 
In our SCFT calculations, we choose a widely used ``red-green-refinement''
method to refine and coarsen the mesh element\,\citep{bank1983}. 
To better demonstrate the performance of the \textit{Log} marking strategy, we take
the Poisson equation in L-shape domain with  Dirichlet boundary condition as
an example  and compare it with $L^2$ marking method in the Appendix section.

\subsection{Surface integral}
\label{subsec:spaceint}

In this subsection, we discuss the integral method on an arbitrary curved
surface.
  We use the symmetrical quadrature rules for triangle element
  \citep{williams2014}.
  Notice that,
  the quadrature points are given in the form of barycentric coordinates. 
  Let $\{\bm\lambda_j\}_{j=0}^{\mathcal J-1}$ and
$\{w_j\}_{j=0}^{\mathcal J-1}$ be the quadrature points and weights, respectively.
For each $\bm\lambda_j$, one can easily obtain a unique point $\bx_j\in\tau_h^p$ and
a unique $\bu_j\in\hat{\tau}$, respectively. We use the following integral formula to calculate the
integral of finte element function $v_h$ defined on $\bM_h^p$:
\begin{align*}
    \int_{\bM_h^p} v_h(\bx) \,\mathrm{d}\bx
    = \sum_{\tau_h^p} \int_{\tau_h^p} v_h(\bx)\,\mathrm{d}\bx &=
    \sum_{\tau_h^p} \int_{\hat{\tau}} v_h(\bx(\bu))|\bm J(\bu)|\,\mathrm{d}\bu \\
  &\approx |\hat{\tau}| \sum_{\tau_h^p} \sum_{j=0}^{\mathcal J-1} w_j
    v_h(\bx_j)|\bm J(\bu_j)|,
\end{align*}
where $\bm J(\bu)$ is the Jacobi matrix of one-to-one mapping
$\bx:\hat{\tau} \mapsto \tau_h^p$ and $|\bm J(\bu)|$ is the corresponding
  determinant.

\section{Contour discretization scheme and integration}
\label{sec:contour}

\subsection{Spectral deferred correction (SDC) method}
\label{subsec:sdc}

In this subsection, we use the SDC method to discretize chain contour variable. The
SDC scheme was originally proposed by Dutt \textit{et al.}\,\citep{Dutt2000sdc} to solve
ordinary differential equations with an appropriate basic method, then used the
residual equation to improve the approximation order of numerical solution.
Unlike the classical deferred correction approach\,\citep{Dutt2000sdc}, the key
idea of the SDC method is to use a spectral quadrature (see
Sec\,\ref{subsec:spint}) to integrate the contour derivative, which can achieve
a high-accuracy numerical solution with a largely reduced number of quadrature
points. In 2019, Ceniceros\,\citep{ceniceros2019efficient} firstly introduced the
SDC method into the SCFT calculations. He chose the second-order
implicit-explicit Runge-Kutta method as the basic solver to discretize the
contour derivative, then used the SDC technique to improve the numerical
precision. We follow a similar idea in this work but use a variable step size
Crank-Nicholson method as the basic solver to compute propagator equations. The
concrete implementation is given as follows.

The Crank-Nicholson method with variable time step size for semi-discrete
propagator equation \eqref{eq:matrix} is 
\begin{equation}
\label{eqn:CN}
\bbM\frac{\bq^{n+1}-\bq^n}{\delta t_n} =
-\frac{1}{2}(\bA+\bF)(\bq^{n+1}+\bq^n),\quad n =0, 1, \dots, N_t-1,
\end{equation}
where $\delta t_n=t_{n+1}-t_{n}$ is the time step size, $t_n$
($n=0,1,2,\dots,N_t-1$) is the Chebyshev node\,\citep{clen1960integral}.
Solving the above equation \eqref{eqn:CN} can obtain the initial numerical solution
$\bq^{[0]}(t)$.
Then we use the SDC scheme to achieve a high-accuracy numerical
solution. We integrate the semi-discrete equation \eqref{eq:matrix}
along the contour variable $t$
\begin{equation}\label{eqn:integral}
	\bbM\bq(t)=\bbM \bq(0)+\int_{0}^{t}\left[(-\bA- \bF) \bq(\tau)\right] d \tau.
\end{equation}
The error between the numerical solution $\bq^{[0]}(t)$ and the exact
semi-discrete solution $\bq(t)$ is 
\begin{align}
	\be^{[0]}(t) = \bq(t) - \bq^{[0]}(t).
\end{align}
The error integration equation
\begin{equation}\label{eqn:error}
\begin{aligned}
    \bbM\be^{[0]}(t)	= & \int_{0}^{t}\left[(-\bA- \bF) \be^{[0]}(\tau) \right] d \tau +
	\boldsymbol{\gamma}^{[0]}(t),
\end{aligned}
\end{equation}
can also be solved by the variable step size CN scheme\,\eqref{eqn:CN} and the
residual equation
\begin{equation}
	\boldsymbol{\gamma}^{[0]}(t) =
	\boldsymbol{M}\boldsymbol{q}(0)+\int_{0}^{t}
	\left[\left(-\boldsymbol{A}
	-\boldsymbol{F}\right)\boldsymbol{q}^{[0]}(\tau)\right]
	d\tau-\boldsymbol{M}\boldsymbol{q}^{[0]}(t),
	\label{eq:residual}
\end{equation}
is computed by the Clenshaw-Curtis quadrature integral method as presented in
Sec.\,\ref{subsec:spint}.
Then the corrected numerical solution is 
\begin{align}
	\bq^{[1]}(t) = \bq^{[0]}(t) + \bm{\be}^{[0]}(t).
\end{align}
Repeating the above process, one can obtain $\bq^{[2]},\dots,\bq^{[J]}$,
$J$ is the pre-determined number of deferred corrections. 
The convergent order of deferred correction solution along the contour variable
is 
\begin{align} 
  \|\bq(t) - \bq^{[J]}(t)\| = O( (\delta t)^{m(J+1)}), 
\end{align} 
where $\delta t = \max\{\delta t_n\}_{n=0}^{N_t-1}$, $m$ is the
order of the basic numerical scheme to solve Eqns.\,\eqref{eq:matrix} and
\eqref{eqn:error}. For the CN scheme, $m=2$. 

\subsection{Spectral integral method along the contour variable} 
\label{subsec:spint} 

In this subsection, we present the Clenshaw-Curtis quadrature
method\,\citep{ceniceros2019efficient} to solve the residual equation
\eqref{eq:residual} and evaluate the density functions \eqref{eq:phiA} and
\eqref{eq:phiB}, which has the spectral accuracy
for smooth integrand function\,\citep{trefethen2008guass}. 
Concretely, assume $g(t)$, $t\in[a,b]$ is a smooth function, the Chebyshev-Gauss
quadrature scheme is 
\begin{equation}
  \begin{aligned}
	\int_{a}^{b} g(t)\,dt &= \frac{b-a}{2}\int_{-1}^{1}g(s)\,ds
    \\
	&\approx
	\left\{\begin{array}{ll}\dfrac{b-a}{2}\left[ a_0 +
            \sum\limits_{k=2 \atop k\, \text{is\ even}}^{{N_t}-2}
        \left(\dfrac{2a_{k}}{1-k^{2}}\right)+\dfrac{a_{N_t}}{1-N_t^{2}}\right],
			& \quad N_t  \ \text{is even},\\
		\\
        \dfrac{b-a}{2}\left[ a_0 + \sum\limits_{k=2 \atop k\, \text{is\ even
    }}^{{N_t}-1}\left(\dfrac{2 a_{k}}{1-k^{2}}\right)\right], & \quad N_t \
    \text{is odd},\\
\end{array}\right.
  \end{aligned}
\end{equation}
where $a_k$, $k=0,1,\cdots,N_t$, is the discrete cosine transform
coefficient of $g$.

\section{Numerical results}
\label{sec:rslt}

In this section, we give several numerical examples to demonstrate the
performance of the proposed methods. In the implementation, we use the linear
($p=1$), quadric ($p=2$), cubic ($p=3$) surface FEM, and adaptive surface FEM,
to discretize the spatial variables, and the SDC method with a one-step
correction to discretize the contour variable. 
The iteration step size $\lambda_{\pm}=2$ in the alternative direction
explicit Euler method in the following computations. The following numerical example
is implemented through FEALPy\,\citep{fealpy2021}.

\subsection{Efficiency}
\label{subsec:efficiency}
In this subsection, we use a parabolic equation with an exact solution 
and SCFT model on a spherical surface to examine the effectiveness of our numerical
methods.

\subsubsection{Efficiency of solving parabolic equation}
\label{subsubsec:effHeat}

As discussed above, the most computationally demanding part of the SCFT simulation is
solving the PDEs for propagators which is a parabolic equation. We take a
parabolic equation with an exact solution as an example to demonstrate the
accuracy of our methods. Consider the parabolic equation defined on a unit sphere
\begin{equation}\label{eq:heat}
	\left\{\begin{array}{ll}{\dfrac{\partial }{\partial t}u(x,y,z,t)=
			\Delta_{\bM} u(x,y,z,t)+f(x,y,z,t)},
			\\
		\vspace{0.2cm}
		{u(x,y,z,0)=\sin (\pi x) \sin (\pi y) \sin (\pi z),} 
	\end{array}\right.
\end{equation}
where $(x,y,z)\in \{(x,y,z)~|~x^2+y^2+z^2=1\},~ t\in[0,T]$, and 
\begin{equation}
\begin{aligned}
    f(x,y,z,t) &= \left(2\pi^2 -1\right)e^{-t} \sin (\pi x) \sin (\pi y) \sin (\pi z) \\
    & + \frac{2 \pi e^{-t} \left[ z \sin (\pi x) \sin (\pi y) \cos (\pi z) + y \sin
	(\pi x)  \sin (\pi z)
    \cos (\pi y) \right]}{x^2 + y^2 + z^2}\\
    & +  \frac{2 \pi e^{-t} \left[ x\sin (\pi y) \sin (\pi z) \cos (\pi x) + xy \sin (\pi z)  \cos (\pi x) \cos (\pi y) \right]}{x^2 + y^2 + z^2}\\
    &+ \frac{2 \pi^2 e^{-t} \left[x z
	\sin (\pi y) \cos (\pi x)
    \cos (\pi z) + y z \sin (\pi x) \cos (\pi y) \cos (\pi z) \right]}{x^2 + y^2 + z^2},
\end{aligned}
\end{equation}
with exact solution $u(x,y,z,t) = \sin (\pi x) \sin (\pi y) \sin(\pi z) e^{-t}$.

\begin{figure}[htbp!]
	\centering
	\setlength{\abovecaptionskip}{0.cm}
	\setlength{\belowcaptionskip}{-0.cm}
	\begin{minipage}[!htbp]{0.4\linewidth}
        \includegraphics[width=5cm]{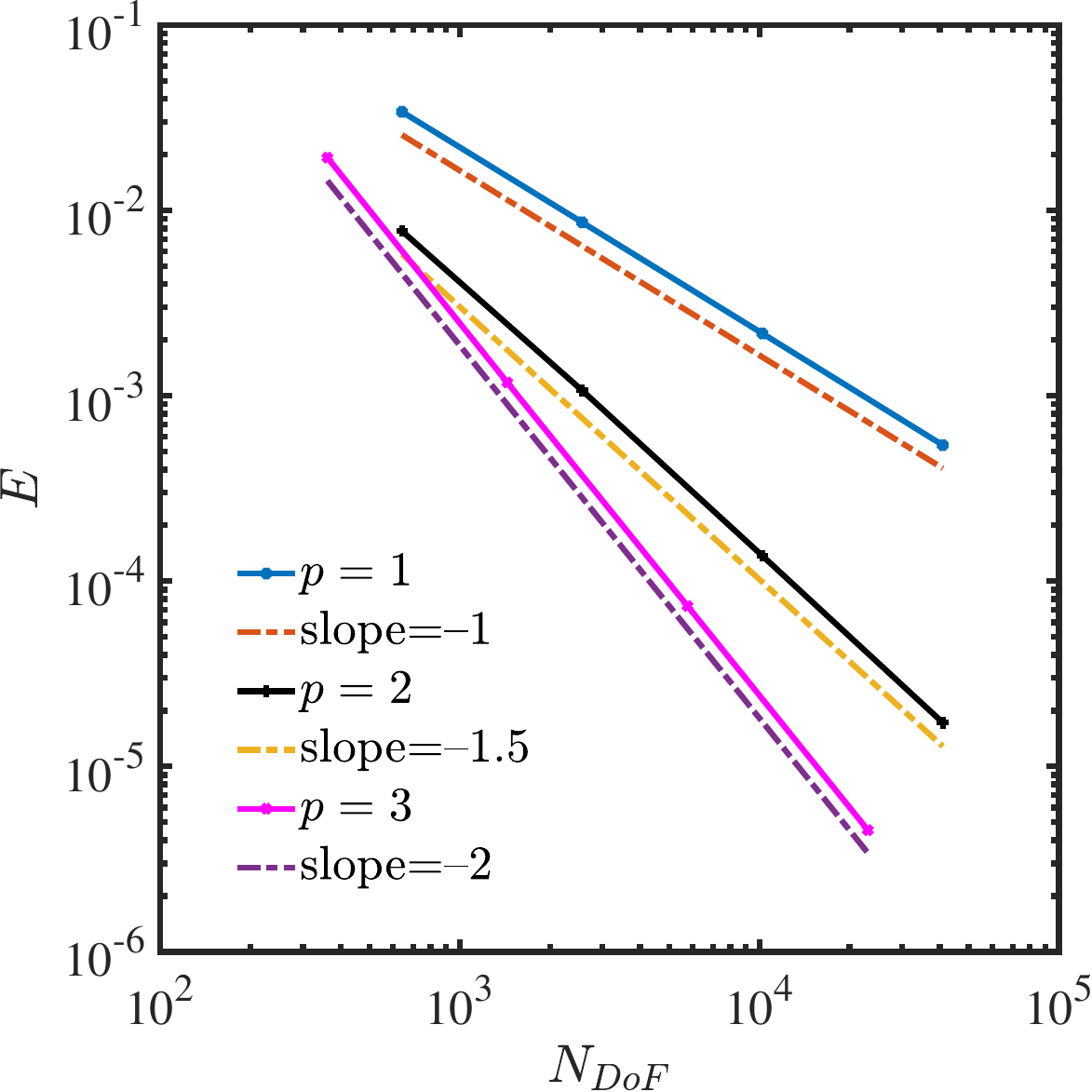}
		\caption*{(a)}
	\end{minipage}
	\hspace{0.6in}
	\begin{minipage}[!htbp]{0.4\linewidth}
        \includegraphics[width=5cm]{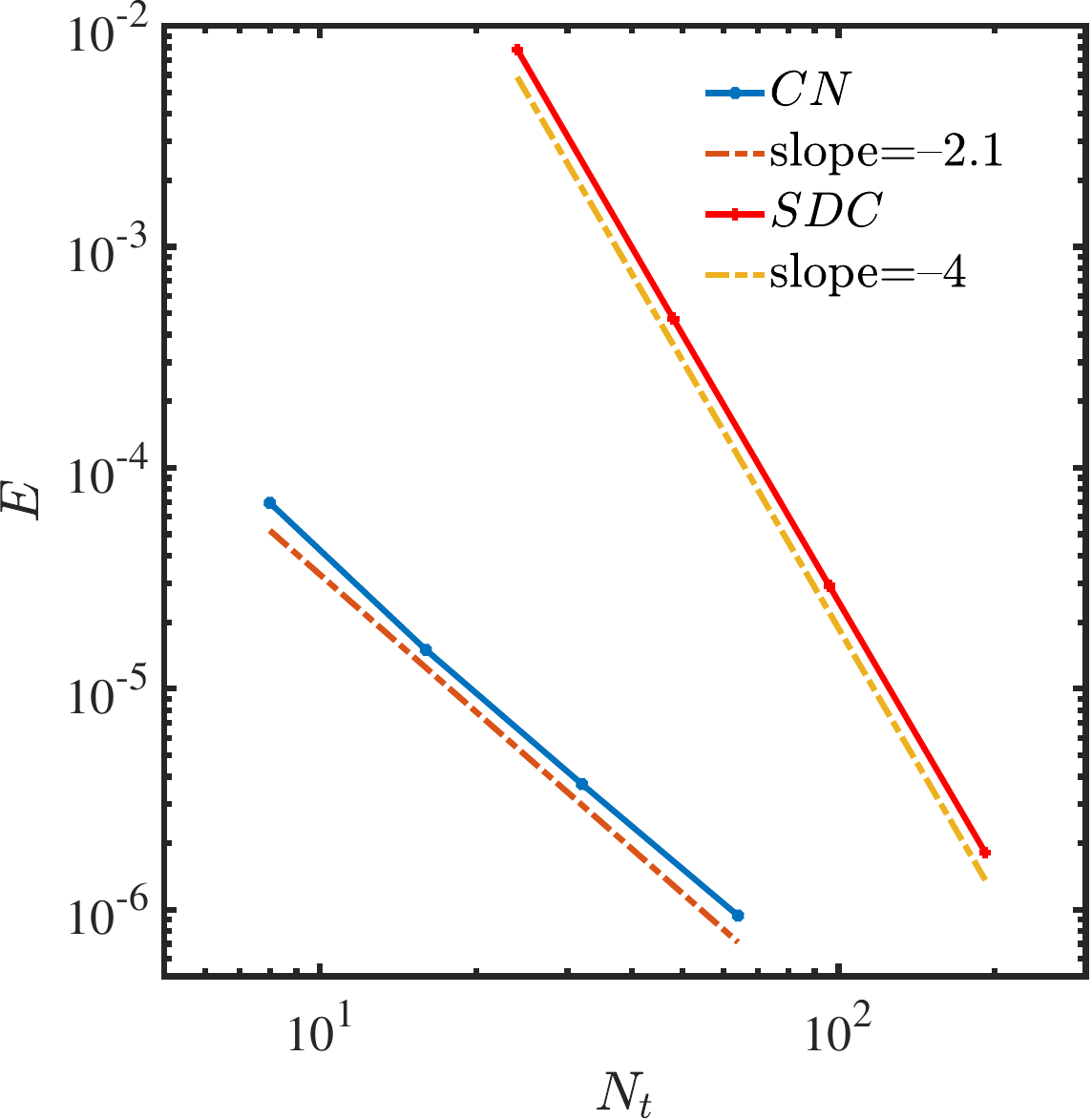}
	\caption*{(b)}
	\end{minipage}
    \caption{\label{fig:Order} Comparison of logarithmic errors against the
	discretization nodes. (a) The spatial error of surface FEMs. $N_{DoF}$ is spatial
	DoFs. (b) The contour error order of CN and SDC
    schemes. $N_t$ is the number of contour points. $E=\|u_I(\cdot, T) -
      u_h(\cdot, T)\|_{\bM^p_h}$, where $u_I(\cdot, T)$ is the value of  exact
      solution $u(\cdot, T)$
    interpolated from $\bM$ to $\bM^p_h$, $T=1$. The error order is the absolute value of the
slope multiply by the dimension of variables, and the
dimensions of space and time are $2$ and $1$, respectively.
}
\end{figure}

We analyze the spatial approximation order of surface FEMs and time
approximation error obtained by CN and SDC schemes. 
Denote $\|\cdot\|_{\bM}$ as the $L^2$ norm on a curved
surface $\bM$. Firstly, we observe the error of surface FEMs. 
We use the CN scheme with $\delta t=1\times 10^{-4}$ to ensure the contour
discretization accuracy.
As shown in Fig.\,\ref{fig:Order}~(a), the error order of linear, quadratic and
cubic surface FEMs are $2$, $3$, and $4$, respectively, which is consistent with
theoretical results.

Secondly, we present the error order of the contour discretization schemes, see
Fig.\,\ref{fig:Order}~(b).  We use the cubic surface FEM as the spatial
discretization scheme to guarantee enough spatial discretization accuracy.
When using the CN method, $5762$ initial spatial degrees of freedom (DoFs) are used.  
When applying the SDC scheme, the initial spatial DoFs are $362$. 
In these experiment, the spatial mesh grid will be refined once as the contour
nodes are doubled. 
From Fig.\,\ref{fig:Order} (b), one can find that the CN method has second-order
precision. SDC method has $m(J+1)$-order precision, $m$ is the order of basic
numerical method (CN), $J$ is the step of correction, here $m=2$, $J=1$.  These
results are also consistent with theoretical results. 

Finally, we confirm the integral accuracy of the numerical solution along with the
contour variable $t$. We use a modified fourth-order integral
method, see (26) in\,\citep{wei2019finite}, and the spectral integral scheme, as
presented in the Sec.\,\ref{subsec:spint}, to integrate $u^{CN}_{h}$ and
$u^{SDC}_h$ from $0$ to $1$ for $t$ to obtain $U_h^{CN}$ and
$U_h^{SDC}$, respectively. $u^{CN}_{h}$ and $u^{SDC}_h$ are the numerical
solutions obtained by CN and SDC methods. The cubic surface FEM with $92162$
DoFs are used to ensure enough spatial precision. 
The integral value of the corresponding exact solution is discretized as $U^e_h$.
The error is defined as 
\begin{align}
  e^{W} = \| U^e_h - U^{W}_h\|_{\bM_h^p},
\end{align}
where $W\in\{CN, SDC\}$. As Tab.\,\ref{tab:intOrder} presents, the SDC method 
achieves the accuracy about $1.89 \times 10^{-7}$ only requiring $16$ contour
points, while the CN scheme needs $128$ nodes to reach the accuracy of $4.65\times
10^{-7}$. It is worth noting that the error value of SDC method is only
reduced to about $2 \times 10^{-7}$ due to the limitation of spatial
discretization precision.
\begin{table}[hbtp!]
	\centering
	\caption{\label{tab:intOrder}
	The contour integral error of the 4-order integral scheme and
	the spectral integral method, $N_t$ is the number of contour points.}
	\begin{tabular}{ccc}
		\hline
		$N_t$&$e^{CN}$&$e^{SDC}$\\
		\hline
        8&5.08e-05&1.02e-06\\
        16&1.32e-05&1.89e-07\\
        32&3.52e-06&2.35e-07\\
        64&1.08e-06&2.53e-07\\
        128&4.65e-07&2.58e-07\\
		\hline
	\end{tabular}
\end{table}

\subsubsection{Efficiency of SCFT calculations}
We apply the high-order adaptive surface FEMs and contour
discretization schemes to solve SCFT equations
to demonstrate the power of our numerical methods.
Due to the complicated SCFT
self-consistent field system, we choose the value of the single chain partition
function $Q$ as a metric to compare the accuracy of different numerical
methods, since it is the integral value of the final result of the propagators.
Here we use a sphere with a radius of $3.56 R_g$ as the calculation surface, a
high-order triangular mesh, see a schematic plot in
Fig.\,\ref{fig:mpmesh}\,(b), to approximate the spherical surface.
We obtain a spotted phase when parameters $\chi N = 25$, $f = 0.2$,
as shown in Fig.\,\ref{fig:estimator}~(b).

\begin{figure}[htbp!]
	\centering
	\setlength{\abovecaptionskip}{0.cm}
	\setlength{\belowcaptionskip}{-0.cm}
	\begin{minipage}[!htbp]{0.4\linewidth}
		\includegraphics[width=5cm]{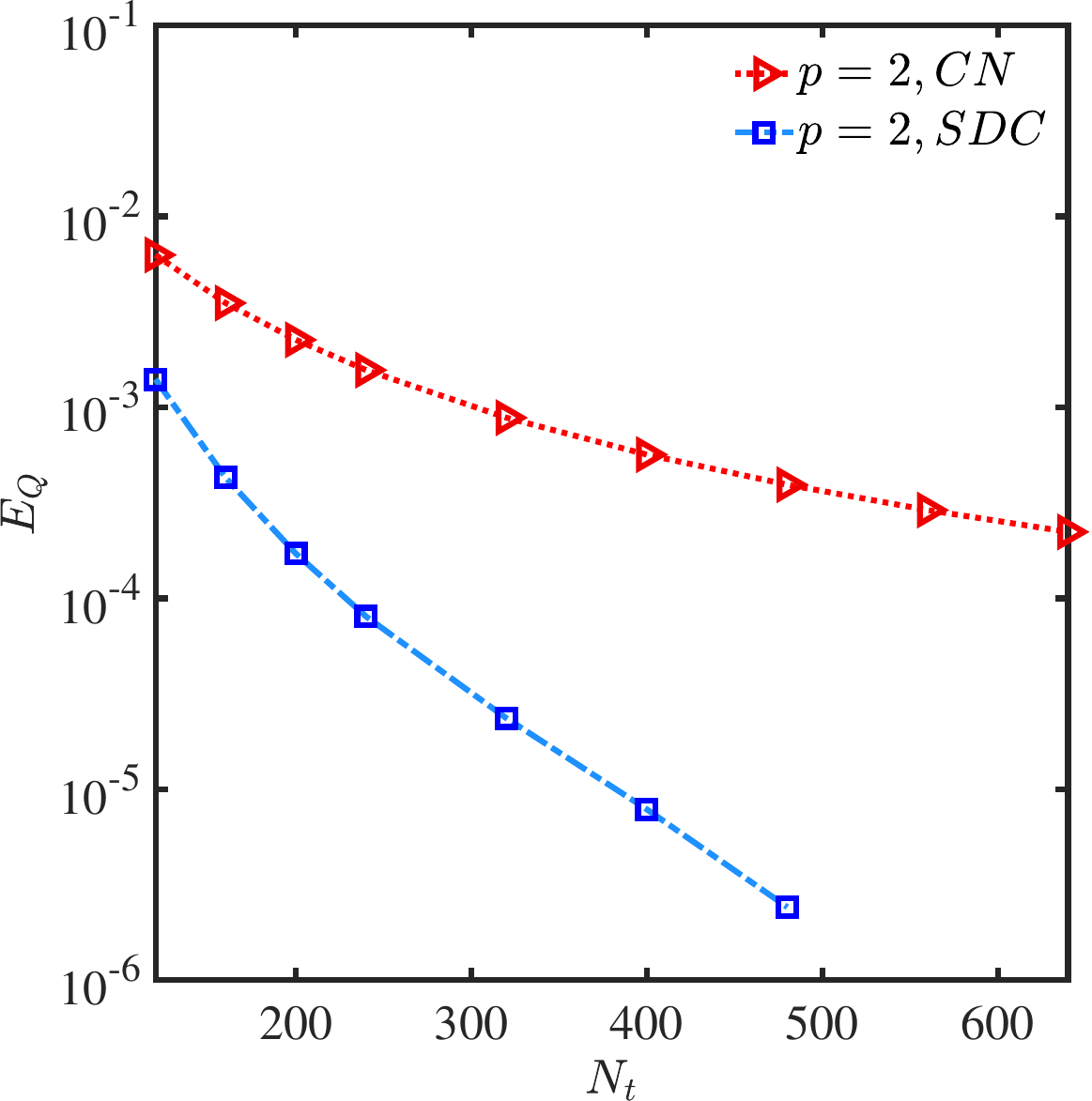}
		\caption*{(a)}
	\end{minipage}
	\hspace{0.6in}
	\begin{minipage}[!htbp]{0.4\linewidth}
		\includegraphics[width=5cm]{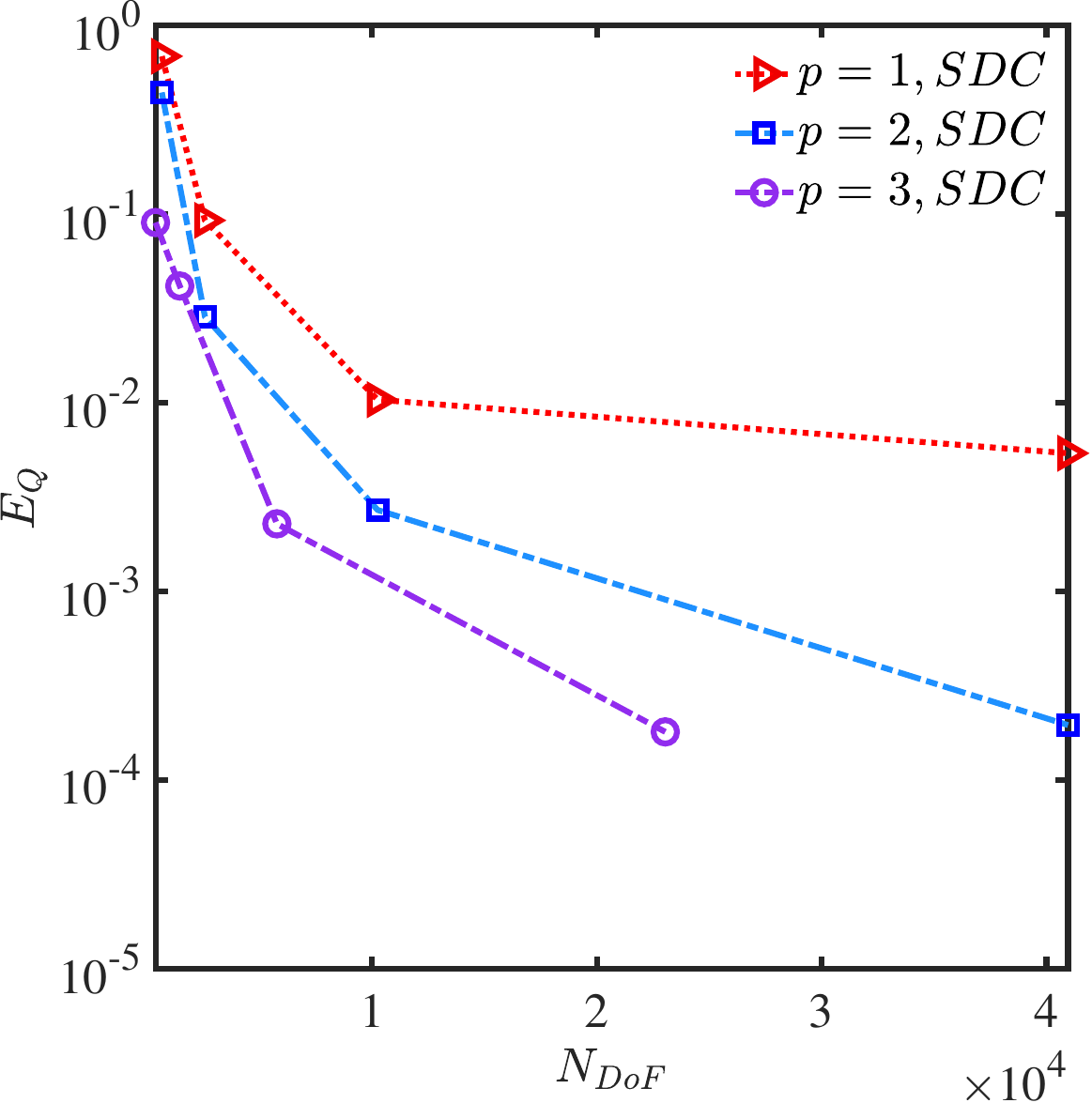}
	\caption*{(b)}
	\end{minipage}
	\caption{\label{fig:scfteff} The convergence of single chain partition function
        $Q$ computed by different numerical schemes. 
	$E_Q$ is the relative error. $N_{DoF}$ is spatial DoFs.
	(a) $Q$ obtained by the CN and SDC scheme as the contour
	point $N_t$ increases, and quadratic surface FEM with $163842$ DoFs
    is employed to discretize the spatial variables.
	(b) $Q$ computed by the linear, quadratic, and cubic surface FEMs with an
	increase of $N_{DoF}$, and the SDC scheme with $200$ points is applied to discretize the contour variable.}
\end{figure}

Firstly, we show the effectiveness of the contour discretization schemes. 
We use quadratic surface FEM with $163842$ DoFs to ensure enough spatial discretization accuracy. 
Certainly, one can choose linear and cubic FEMs, as long as the spatial
discretization is accurate enough.
In Fig.\,\ref{fig:scfteff} (a), $Q_{ref}$ is numerically
calculated by the
SDC method with  $560$ contour points. $E_Q =(Q-Q_{ref})/Q_{ref}$ is the relative
error. Fig.\,\ref{fig:scfteff} (a) shows that the SDC method converges much faster than the CN method.

Secondly, we investigate the accuracy of $p$-order ($p=1,2,3$) surface FEMs in the
SCFT simulation. 
We use the SDC scheme with $200$ discretization points to guarantee enough accuracy in contour direction.
In Fig.\,\ref{fig:scfteff} (b), the reference value
$Q_{ref}$ is numerically obtained by the cubic surface FEM with $92162$ DoFs.
$E_Q =(Q-Q_{ref})/Q_{ref}$ is
the relative error. One can observe that, compared with the linear surface
FEM, the quadric and cubic surface FEMs have faster convergent rates.

\subsection{Adaptive surface FEM}
\label{subsec:rltsAdaptiv}

In this subsection, we illustrate the performance of adaptive mesh method from two
parts: the computational cost to achieve the same precise level
and the application to strongly segregated systems compared with the uniform
mesh approach.  In this subsection,
we choose quadratic surface FEM in the adaptive mesh method.
Certainly, one can apply arbitrary order surface FEMs in this adaptive mesh approach.
The SDC scheme with $200$ contour points is used to guarantee
the contour discretization accuracy.

\begin{figure}[htbp!]
	\centering
	\includegraphics[width=0.8\linewidth]{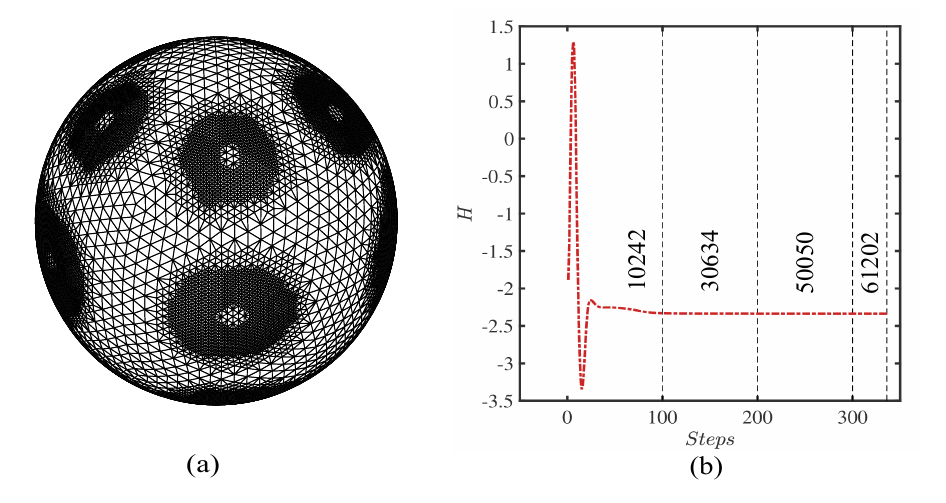}
	\caption{\label{fig:afem:eff}
        (a) The final adaptive mesh. (b) The change of Hamiltonian $H$ during the
		adaptive procedure.
	The numbers between two dotted lines represent
	spatial DoFs in the adaptive process.}
\end{figure}

Firstly, we discuss the effectiveness of adaptive mesh method through an example in
which the parameter $\chi N = 25$, $f=0.2$, and a computational domain is a sphere with a radius of $3.56 R_g$. The triangular uniform mesh
with $10242$ DoFs is used as the initial mesh of the system.
The adaptive process begins when the SCFT iteration reaches the maximum step
$500$ or the reference estimator satisfies 
$$
e_{ref} = \sigma(e_{\tau})/(\max(e_{\tau})-\min(e_{\tau})) < 0.1,
$$
where $\sigma(e_{\tau})$ is the standard deviation of $e_{\tau}$, as defined in 
Eqn.\,\eqref{eq:afem:error}. The stop condition of the adaptive process
is the Hamiltonian discrepancy less than $1.0\times 10^{-6}$.
The final convergent structure is the spotted phase, as shown in
Fig.\,\ref{fig:estimator}~(b).
Fig.\,\ref{fig:afem:eff}~(a) presents the final adaptive mesh with $61202$ DoFs.
From these results, one can find that the adaptive approach can
detect the sharp interface and damped internal structure to assign the mesh
distribution reasonably.
The change of Hamiltonian $H$ during the adaptive process is given in 
Fig.\,\ref{fig:afem:eff}~(b).

\begin{figure}[!hbpt]
	\centering
	\setlength{\abovecaptionskip}{0.cm}
	\setlength{\belowcaptionskip}{-0.cm}
	\begin{minipage}[!htbp]{0.4\linewidth}
		\includegraphics[width=5cm]{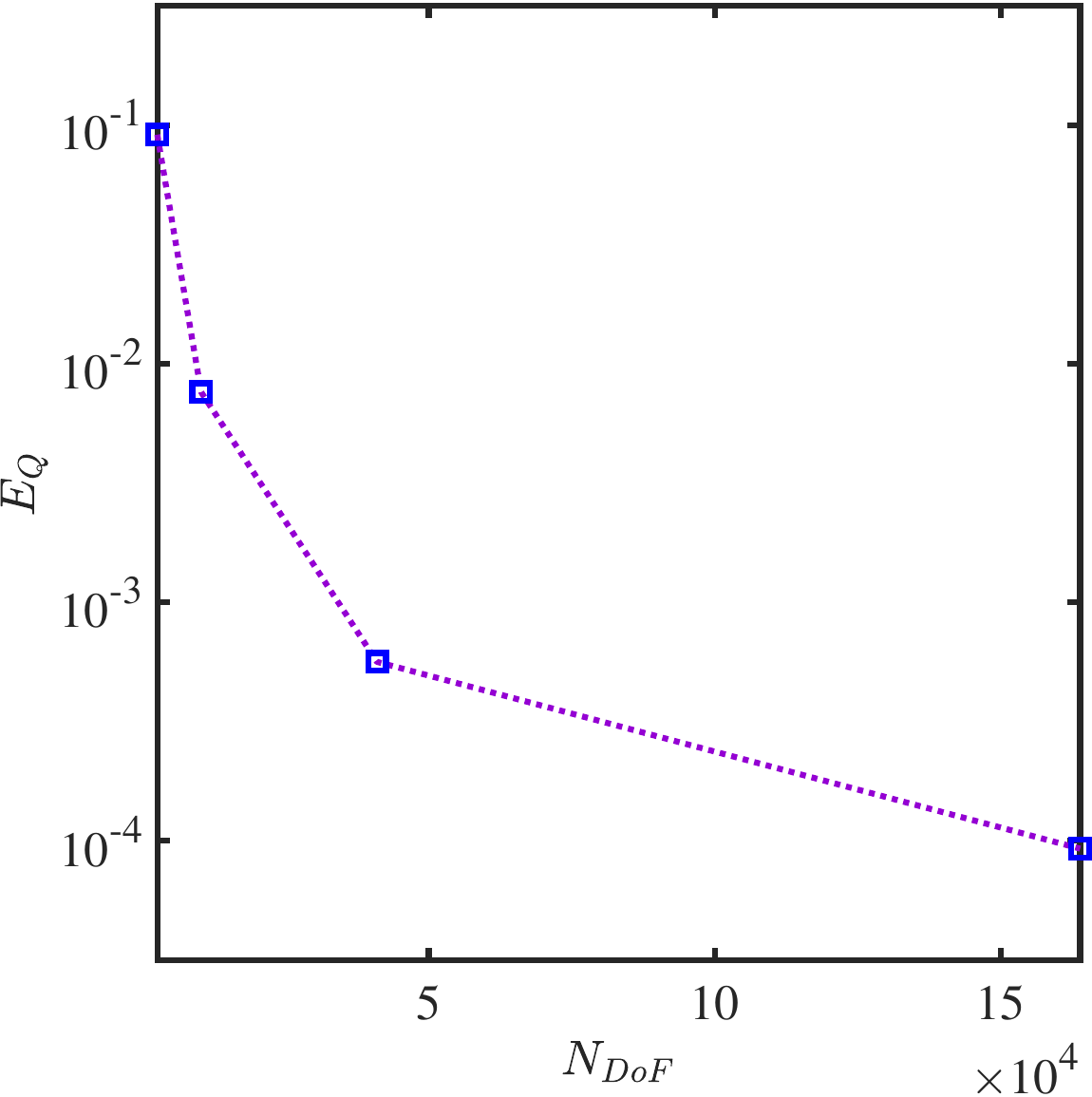}
		\caption*{(a)}
	\end{minipage}
	\hspace{0.6in}
	\begin{minipage}[!htbp]{0.4\linewidth}
		\includegraphics[width=5cm]{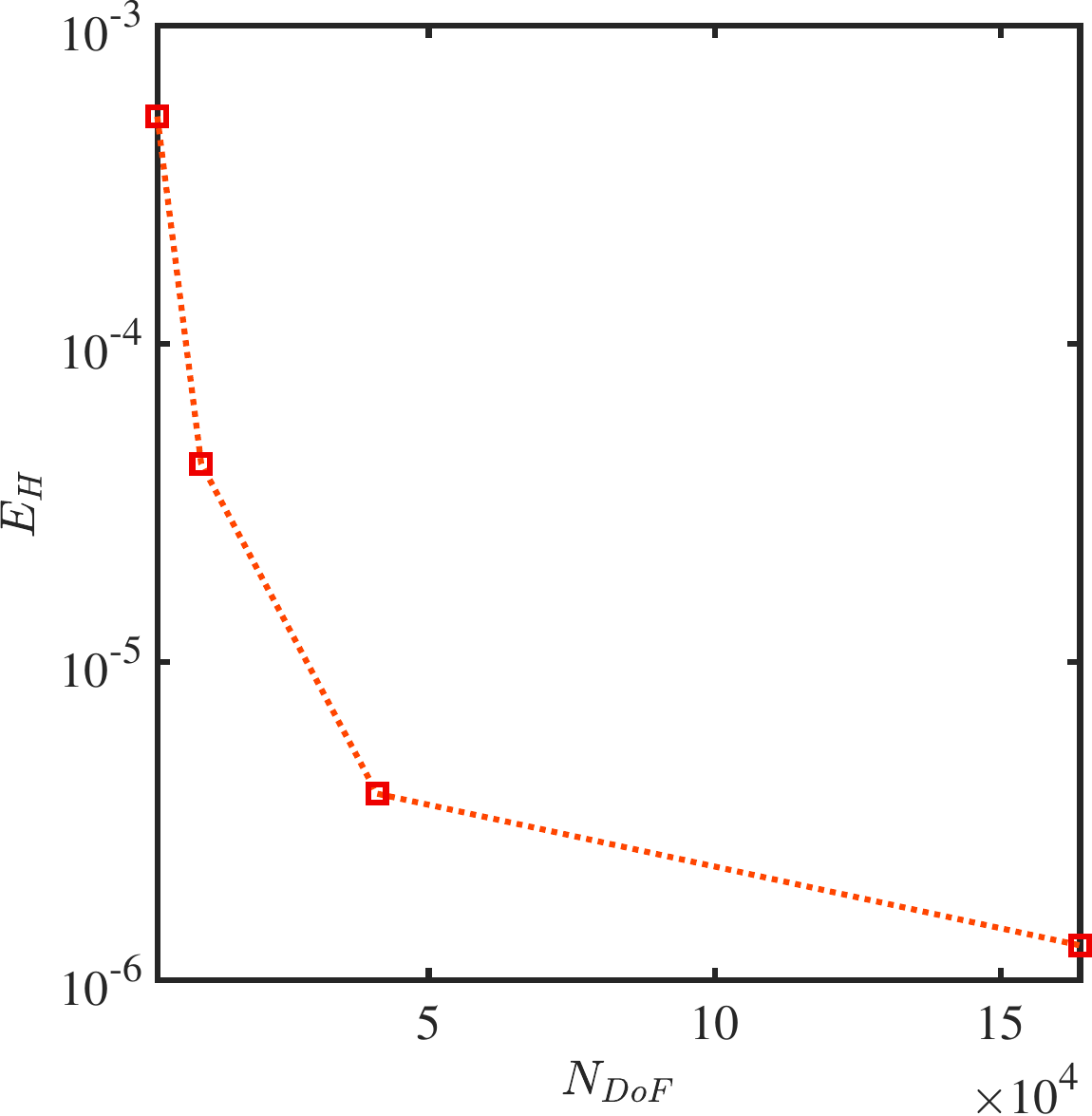}
	\caption*{(b)}
	\end{minipage}
	\caption{\label{fig:compare} 
	A comparison between the adaptive mesh method and the uniform mesh
	approach when $\chi N =25, f=0.2$.  The relative error (a) $E_Q$ and (b)
    $E_H$.}
\end{figure}

Secondly, we further compare the values of $Q$ and $H$ when using the adaptive mesh
and uniform mesh approaches. Fig.\,\ref{fig:compare} shows the convergence of $Q$ and $H$
by the uniform mesh method. $E_Q =(Q-Q_{ref})/Q_{ref}$ and $E_H =
(H-H_{ref})/H_{ref}$ are the relative errors. The reference values $Q_{ref}$ and
$H_{ref}$ are obtained by the adaptive mesh method with $61202$ DoFs. 
Tab.\,\ref{tab:afem-uniform} presents the finally converged values of $Q$ and $H$
through the adaptive mesh and uniform mesh methods. One can observe
that $Q$ and $H$ obtained by the uniform mesh method converge to the adaptive mesh approach
from these numerical results. The uniform mesh method needs $163842$ DoFs to achieve
the same error level, which is $2.68$ times than the adaptive mesh approach does.

\begin{table}[htbp!]
	\centering
	\caption{\label{tab:afem-uniform} 
    The converged values of $Q$ and $H$ otained by the adaptive and uniform mesh
	methods when $\chi N = 25, f = 0.2$.}
	\begin{tabular}{cccc}                                                
		\hline
        Method &$N_{DoF}$ & $Q$ & $H$\\
		 \hline
		Adaptive &  61202   & 3.2329e+02 & -2.336168\\
		Uniform  &  163842 & 3.2326e+02 & -2.336165\\
		\hline
\end{tabular}
\end{table}

\begin{figure}[!htbp]
	\centering
	\includegraphics[width=1.0\linewidth]{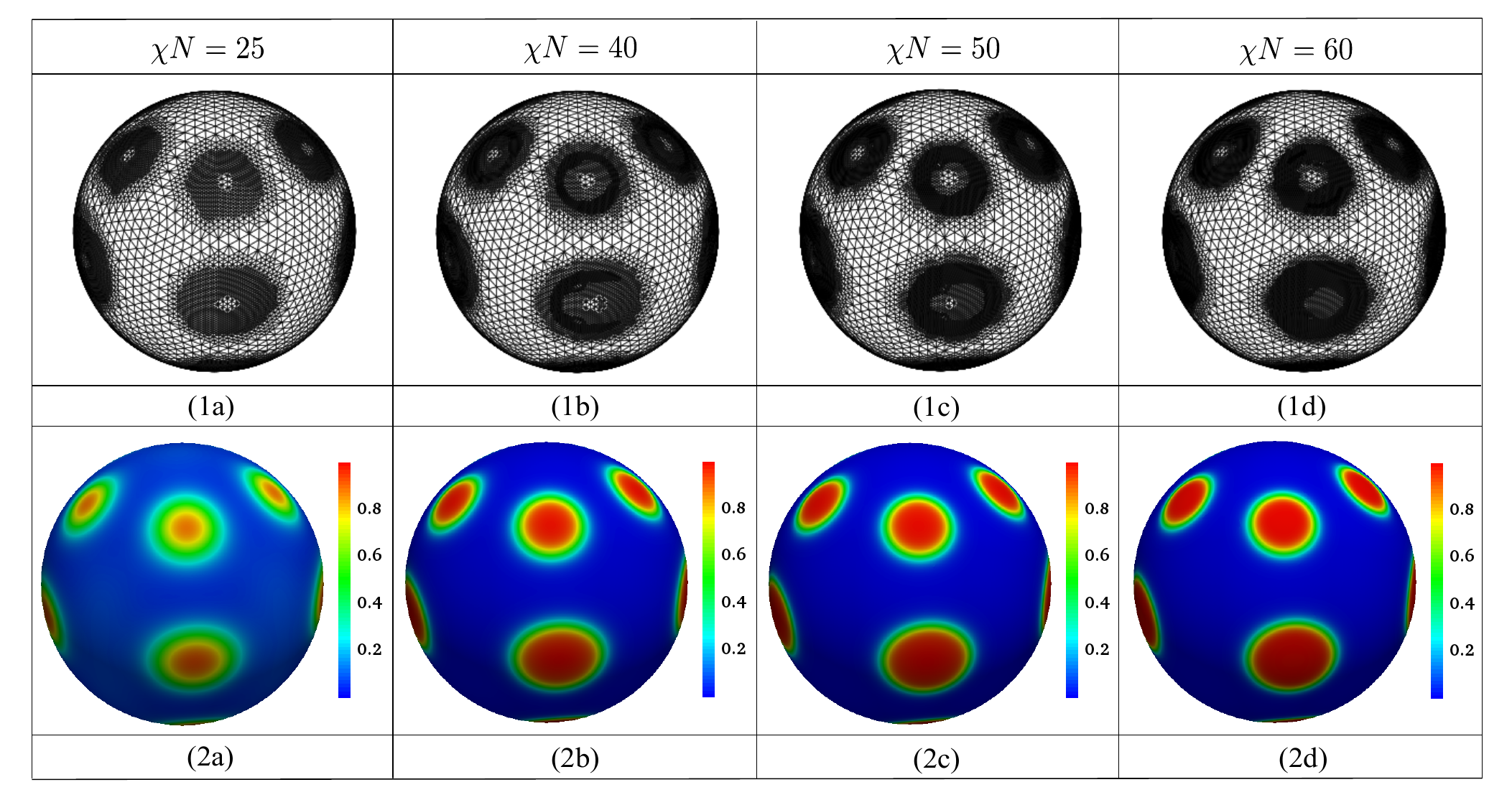}
	\caption{\label{fig:afem:spere} The adaptive meshes and equilibrium structures of
$f=0.2$ and different interaction parameter $\chi N$.  }
\end{figure}

Thirdly, we apply the adaptive mesh method to the strongly segregated systems on a
spherical surface with a radius of $3.56 R_g$.  Tab.\,\ref{tab:afem-adaptive}
gives numerical results of $\chi N$ from $25$ to $60$ and $f=0.2$. In
Tab.\,\ref{tab:afem-adaptive}, $N^{Adap}_{DoF}$ is the spatial DoFs of the adaptive
mesh scheme. $h_{min}$ is the minimum size of the adaptive mesh. $N^{Uni}_{DoF}$
is the estimated DoFs of the uniform mesh method by making the mesh size equal
$h_{min}$.  In relatively strong segregation regime, its initial values come from
the converged results of relatively weak segregated systems. For example, the
initial values of $\chi N =30$ system come from the converged results of $\chi N
=25$ system by the interpolation approach.  
\begin{table}[!hbpt]                                                                
	\centering
	\caption{\label{tab:afem-adaptive}
	Numerical results by the adaptive surface FEM for strongly segregated
    systems with different $\chi N$ and $f=0.2$.
    $N^{Adap}_{DoF}$ is spatial DoFs of the adaptive mesh method. $h_{min}$ is the minimum size
	of the adaptive mesh element. $N^{Uni}_{DoF}$ is the estimated spatial DoFs of the
	uniform mesh scheme with mesh size $h_{min}$.} 
	\begin{tabular}{cccccc}                                                
		\hline
        $\chi N$ & SCFT iterations & $h_{min}$ & $N^{Adap}_{DoF}$ &
		$N^{Uni}_{DoF}$ &
        $N^{Adap}_{DoF}/N^{Uni}_{DoF}(\%)$  \\
        \hline
        25 & 336 & 1.72e-02 & 61202  & 163842 & 37.3\% \\
        30 & 40  & 1.72e-02 & 68130  & 163842 & 41.5\% \\
        35 & 70  & 8.64e-03 & 79146  & 655362 & 12.1\% \\
        40 & 69  & 8.64e-03 & 113554 & 655362 & 17.3\% \\
        45 & 72  & 8.64e-03 & 153634 & 655362 & 23.4\% \\
        50 & 76  & 8.64e-03 & 178058 & 655362 & 27.1\% \\
        55 & 79  & 8.64e-03 & 198290 & 655362 & 30.2\% \\
        60 & 49  & 4.32e-03 & 219634 & 2621442 & 8.37\% \\
		\hline
\end{tabular}
\end{table}
Correspondingly, Fig.\,\ref{fig:afem:spere} shows the distribution of converged
adaptive meshes and equilibrium structures for different $\chi N$. 
It is obvious that the interface narrows as $\chi N$ increases.
From these results, one can find that spatial DoFs and iteration steps of the adaptive mesh method increase mildly as $\chi N$ increases.
The adaptive mesh approach can significantly save the computational amount compared with the uniform mesh approach, even up to 
$91.63\%$ when $\chi N =60$, $f=0.2$. 
This is attributed to the adaptive mesh approach that can efficiently arrange
two-scale mesh grids for the sharp interface and damped internal morphology
by locally refining and coarsening mesh grids.

\subsection{Self-assembled structures on general curved surfaces}
\label{subsec:generalSurfs}

\begin{figure}[htbp!]
	\centering
	\includegraphics[width=1.0\linewidth]{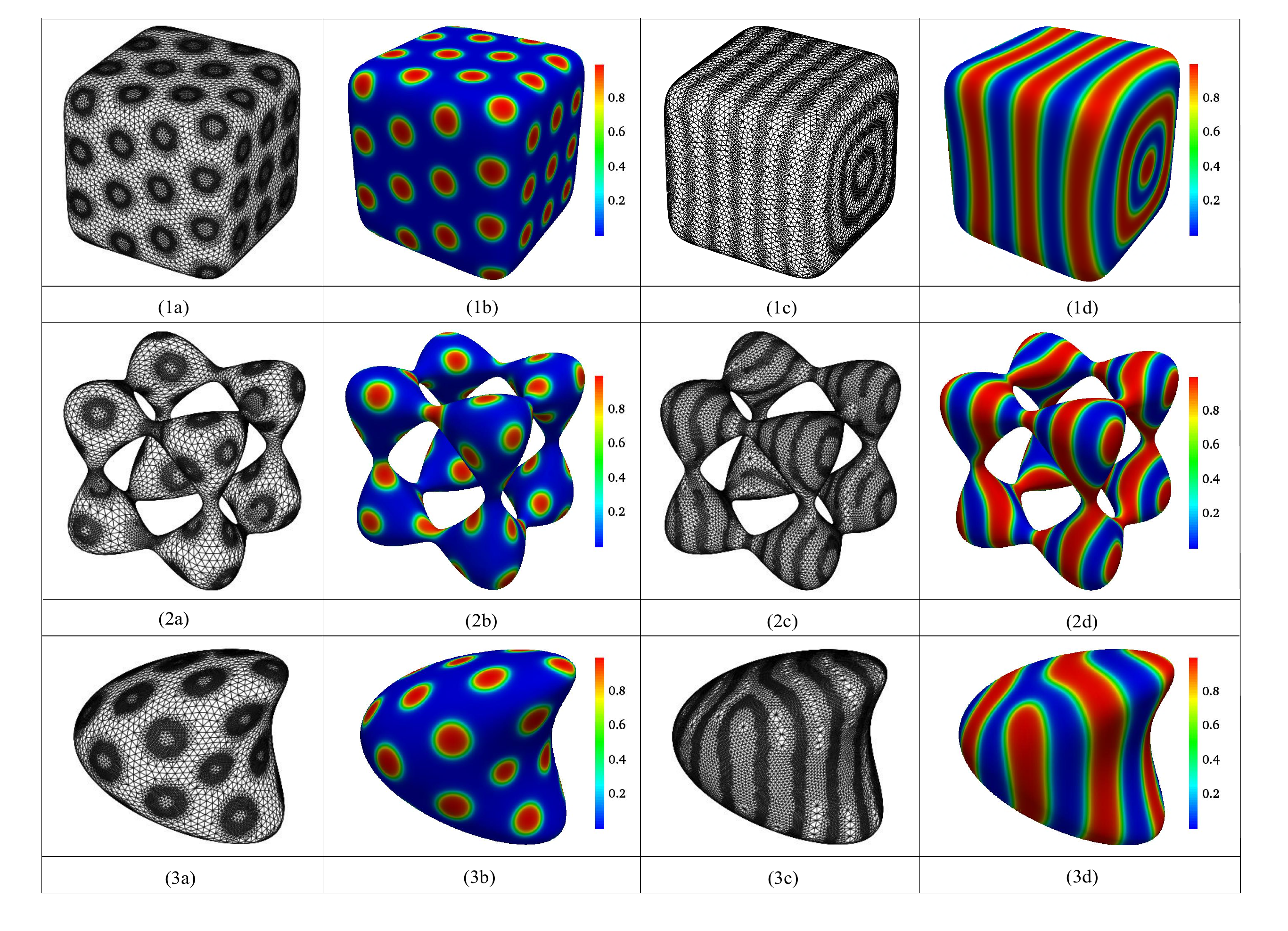}
	\caption{\label{fig:afem:curve} 
	Applying the adaptive mesh method for solving SCFT on three closed curved surfaces 
	(1). Squared surface; (2). Quartics surface; (3). Heart surface. The second and
fourth columns present the spotted phase and strip phase, respectively. The first and
third columns show the corresponding adaptive meshes. Model parameters are
(1b) $[\chi N =45, f=0.2]$, (1d) $[\chi N =20, f=0.5]$, (2b) $[\chi N =35,
f=0.2]$,
(2d) $[\chi N =30, f=0.5]$, (3b) $[\chi N =45, f=0.2]$, (3d) $[\chi N =25,
f=0.5]$. The calculated data are given in the
Tab.\,\ref{tab:data}.}
\end{figure}

\begin{figure}[htbp!]
	\centering
	\includegraphics[width=1.0\linewidth]{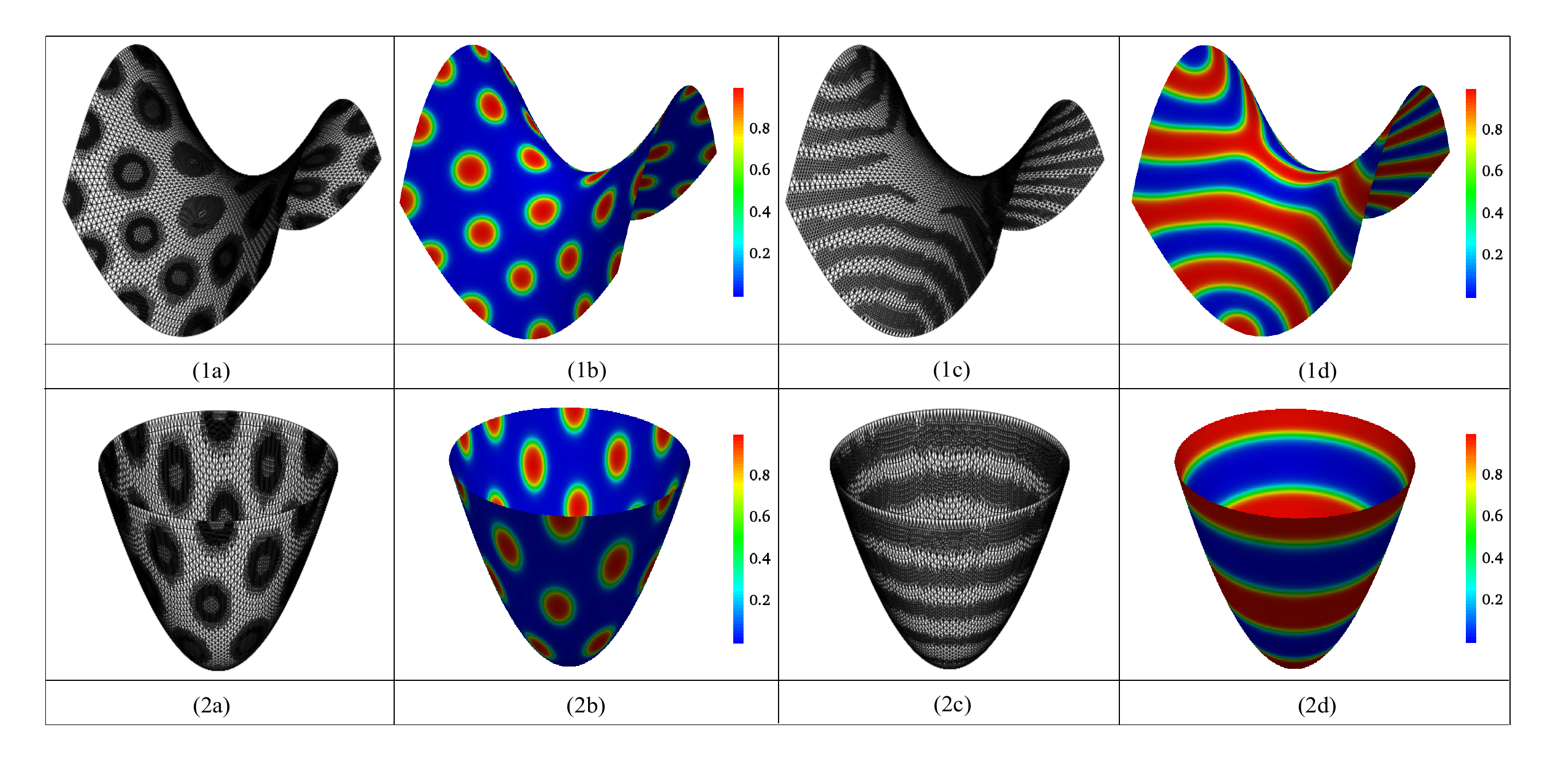}
	\caption{\label{fig:afem:curve1} 
	Applying the adaptive mesh method to solving SCFT on two open curved surfaces 
        (1). Saddle surface and (2). Parabolic surface for relative strongly
        segregated systems. The second and fourth columns present the spotted
        phase and strip phase, respectively. The first and third columns show
        the corresponding adaptive meshes. Model parameters are
(1b) $[\chi N =40, f=0.2]$, (1d) $[\chi N =20, f=0.5]$, (2b) $[\chi N =40,
f=0.2]$,
(2d) $[\chi N =35, f=0.5]$. The calculated data are shown in the
Tab.\,\ref{tab:data}.}
\end{figure}

This subsection presents the application of the adaptive mesh approach on
five curved surfaces for different $\chi N$ and $f$. 
Fig.\,\ref{fig:afem:curve} shows the converged results of adaptive meshes and
morphologies on three closed surfaces, and Fig.\,\ref{fig:afem:curve1} provides
the corresponding results on two open surfaces. The homogeneous Neumann boundary
condition is used for open surfaces. 
Tab.\,\ref{tab:data} presents the corresponding calculated data. From
these results, one can find that the adaptive mesh method is efficient and saves
spatial DoFs both on closed and open surfaces. For example, 
the adaptive mesh method can save the computational amount up
to $89\%$ for computing the spotted phase on the parabolic surface
when $\chi N =40$, $f=0.2$.

\begin{table}[H]                                                                
	\centering
    \caption{\label{tab:data} The computational data on five different curved
        surfaces as shown in Figs.\,\ref{fig:afem:curve} and
		\ref{fig:afem:curve1}.  $N^{Adap}_{DoF}$ and $h_{min}$ are spatial
		DoFs and the minimum element size of the adaptive mesh approach,
		respectively. $N^{Uni}_{DoF}$ is the estimated spatial DoFs of the
		uniform mesh method with mesh size $h_{min}$.}
\begin{tabular}{ccccc}                                                
		\hline
        Mesh & $h_{min}$ & $N^{Adap}_{DoF}$ & $N^{Uni}_{DoF}$ &
        $N^{Adap}_{DoF}/N^{Uni}_{DoF}(\%)$ \\
		\hline
		Fig.\,\ref{fig:afem:curve}(1)(a)  & 7.756e-03 & 444114 & 1483060 &  30\%\\
		Fig.\,\ref{fig:afem:curve}(1)(c)  & 2.458e-02 & 68414  & 93528   &  73\%\\
        Fig.\,\ref{fig:afem:curve}(2)(a)  & 1.129e-02 & 103876 & 290272  &  35\%\\
		Fig.\,\ref{fig:afem:curve}(2)(c)  & 1.039e-02 & 228576 & 356290  &  64\%\\
		Fig.\,\ref{fig:afem:curve}(3)(a)  & 7.799e-03 & 131370 & 689928  &  19\%\\
		Fig.\,\ref{fig:afem:curve}(3)(c)  & 1.505e-02 & 85238  & 172488  &  49\%\\
		Fig.\,\ref{fig:afem:curve1}(1)(a) & 3.907e-03 & 179929 & 1052676 &  17\%\\
		Fig.\,\ref{fig:afem:curve1}(1)(c) & 7.815e-03 & 47685  & 264196  &  18\%\\
		Fig.\,\ref{fig:afem:curve1}(2)(a) & 3.272e-03 & 107919 & 1007300 &  11\%\\
		Fig.\,\ref{fig:afem:curve1}(2)(c) & 3.272e-03 & 267297 & 1007300 &  26\%\\
		\hline
\end{tabular}
\end{table}

\section{Conclusion}
In this paper, we proposed an adaptive high-order surface FEM for solving SCFT
equations on general curved surfaces.  The high-order surface FEM was obtained
through two aspects: high-order function space approximation and high-order
geometrical surface approximation.  To further reduce spatial DoFs, an efficient
adaptive mesh method equipped with a new \textit{Log} marking strategy was
presented, which can make full use of the information of obtaining numerical
results to refine or coarsen mesh.  To improve the approximation order of the
contour derivative, the SDC method was also used to address the contour variable.
The resulting method can achieve high accuracy with fewer spatial and contour
nodes and is suitable for solving strongly segregated systems.  When computing
strongly segregated systems, the superiority of the adaptive mesh method is even
more pronounced, which can save spatial DoFs up to $91.63\%$ at most compared
with the uniform mesh approach.  In the future, we will apply these numerical
methods to other polymer systems, for example, rigid chain systems.

\bibliographystyle{siamplain}

\section*{Appendix: A comparison between \textit{Log} and $L^2$ marking
strategies}

In this Appendix section, we take the Poisson equation \eqref{eq:poisson} as an
example to compare the \textit{Log} with the traditional $L^2$ marking
strategies.
\begin{equation}
\begin{aligned}
  -\Delta u =  0,  & \text{ in } \Omega,\\
  u =  g, & \text{ on } \partial \Omega,
\end{aligned}
\label{eq:poisson}
\end{equation}
where $\Omega = (-1,1)^2\backslash (0,1)\times(-1,0)$. Here we take the exact
solution as $u_e=r^{2/3} \sin(2\theta/3)$, where $r, \theta$ are polar
coordinates. Then the Dirichlet boundary condition $g$ can be determined by $u_e$. 

We use the linear FEM to solve the equation and obtain a numerical solution $u_h$.
Denote the error estimator on the element $\tau$ as
$$
e_{\tau} = \| R_hu_h-\nabla u_h\|_{\tau},
$$
where $R_h u_h$ is the simple average operator
$$
R_h u_h(\bx_i) := \frac{1}{n} \sum_{\bx_i\in \tau, \tau\in \bK} \nabla
u_h(\bx_i)|_{\tau}.
$$
$n$ is the number of element $\tau$ with $\bx_i$ as a vertex, $\bK$ denotes the
set of elements in $\Omega$ with common vertex $\bx_i$.

The \textit{Log} marking strategy has been presented in Sec.\,\ref{subsec:afem}.
In this comparison, we set $\theta=1$ in \eqref{eq:log} in the
\textit{Log} marking strategy.
Before we go ahead, it is necessary to introduce the $L^2$ marking strategy.
Let $\calA$ be a set consisting of all elements which are required to be
refined. The initial $\calA$ is empty. In each adaptive iteration, $L^2$ marking
strategy sorts element estimators $e_{\tau}$ in a descending order, and adds
elements into $\calA$ until
$$
\sum_{\tau \in
\calA}  e_{\tau}^2 \geq \beta\sum_{\tau \in \Omega} e_{\tau}^2, 
$$
where $\beta \in (0,1)$. In the following calculation we set $\beta = 0.2$. We
refine the mesh elements in $\calA$ with ``red-green-refinement'' approach\,\citep{bank1983}. 

\begin{figure}[!hptb]
	\centering
	\setlength{\belowcaptionskip}{-0.cm}
	\begin{minipage}[!htbp]{0.4\linewidth}
        \centering
        \includegraphics[width=4.6cm]{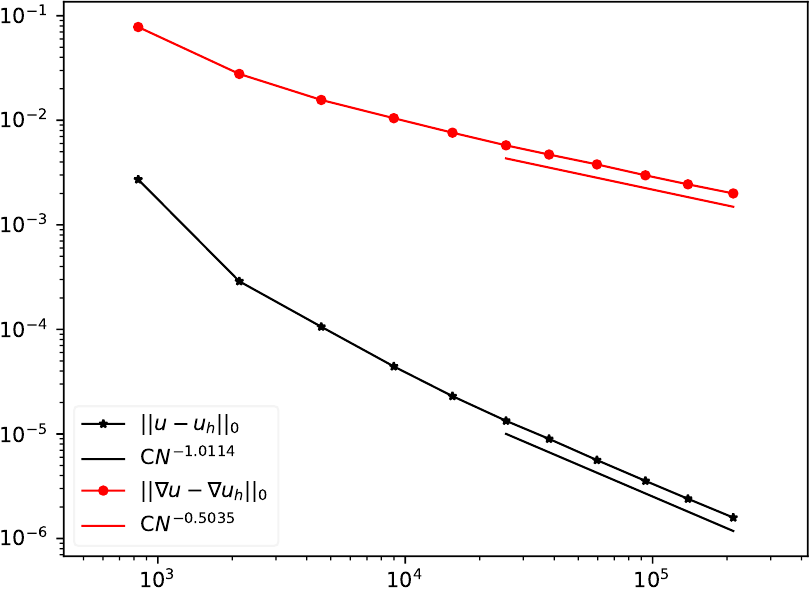}
		\caption*{(a)}
	\end{minipage}
	\hspace{0.6in}
	\begin{minipage}[!htbp]{0.4\linewidth}
        \centering
        \includegraphics[width=4.6cm]{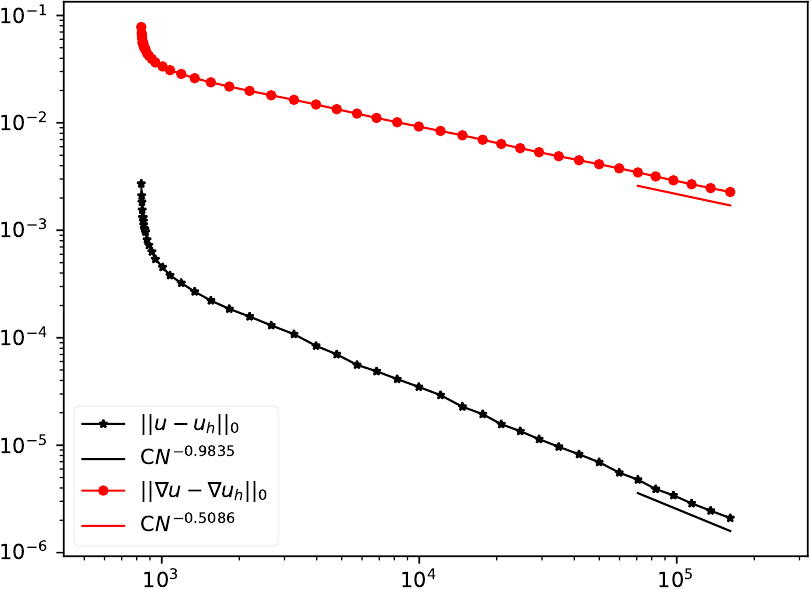}
	\caption*{(b)}
	\end{minipage}
    \caption{\label{fig:markerror}
     The comparison of errors with (a) \textit{Log} and (b) $L^2$ marking
	 strategies. The $x$-axis is $N_{DoF}$.}
\end{figure}

We perform the adaptive method to solve \eqref{eq:poisson} through
\textit{Log} and $L^2$ marking strategies, respectively. 
Numerical results are shown in Tab.\,\ref{tab:lerror} and Fig.\,\ref{fig:markerror}.
We use $\|u_e-u_h\|_0$ with $L^2$-norm to measure the numerical error between
exact solution $u_e$ and numerical solution $u_h$.
When the numerical error reaches an level of about $2.4\times 10^{-6}$, 
two methods have almost same $N_{DoF}$. However, $L^2$ marking strategy requires $46$ 
adaptive iterations, while \textit{Log} strategy only needs $10$ adaptive iterations.
Correspondingly, the \textit{Log} strategy saves $46.2\%$ CPU times.
It is demonstrated that the \textit{Log} strategy can save much CPU time in the adaptive
method.
\begin{table}[H]
	\centering
     \caption{ \label{tab:lerror}
       The relative data for the $L^2$  and \textit{Log}
   marking strategies at the error level $2.4\times 10^{-6}$.}
	\begin{tabular}{ccccc}
	\hline
    Strategy   & $N_{DoF}$ & $\|u_e-u_h\|_{0}$  & Iteration & CPU\\
	\hline
    $L^2$      &   135451   &  2.45e-06      &  46       & 108s \\
    \textit{Log} &  138940   &  2.39e-06      &  10       & 58s \\
	\hline
	\end{tabular}
\end{table}

\begin{figure}[!hptb]
        \includegraphics[width=12.3cm]{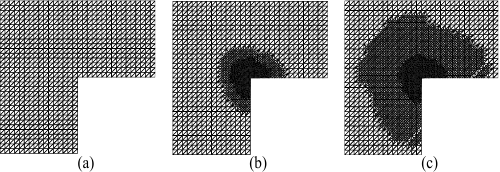}
   \caption{\label{fig:mmesh}
     (a) Initial mesh. Adaptive meshes
     when the error level is $1.0\times 10^{-4}$: (b) the $2^{nd}$
   iteration with \textit{Log} marking strategy; (c) the $26^{th}$
      iteration with $L^2$ marking strategy. 
   }
\end{figure}

Fig.\,\ref{fig:mmesh} shows the adaptive meshes obtained by two strategies when
the error level is $1.0\times 10^{-4}$.  As one can see, to achieve the same
 error level, the $Log$ marking strategy only needs $2$ iterations, while the
 second method needs $26$ iterations. At the same time, the refined meshes of the
 $Log$ marking strategy are concentrated near the singular points. In contrast, 
 the refined meshes of the $L^2$ marking strategy spread to non-singular domain
 which does not need to be refined.
 These adaptive processes sufficiently demonstrate the efficiency of the new
 proposed \textit{Log} marking strategy.

\end{document}